\documentclass{amsart}
\usepackage{amsmath,amssymb,amsthm}

\usepackage[all]{xy}
\usepackage{graphicx}
\usepackage{hyperref}
\usepackage{mathdots}
\usepackage{enumerate}
\usepackage{mathdots}
\usepackage{pstricks}

\theoremstyle{plain}
\newtheorem{The}{Theorem}[section]
\newtheorem{Prop}[The]{Proposition}
\newtheorem{Cor}[The]{Corollary}
\newtheorem{Lem}[The]{Lemma}

\theoremstyle{definition}
\newtheorem{Def}[The]{Definition}
\newtheorem{Ex}[The]{Example}
\newtheorem{Not}[The]{Notation}
\newtheorem{Nots}[The]{Notations}
\newtheorem{Defs}[The]{Definitions}
\newtheorem{Hyp}[The]{Assumption}

\theoremstyle{remark}
\newtheorem{Rem}[The]{Remark}

\newcommand{\C}{\mathbf{C}}

\newcommand{\Z}{\mathbf{Z}}
\newcommand{\N}{\mathbf{N}}

\newcommand{\PP}{\mathbf{P}}
\newcommand{\bi}{\mathbf{i}}
\newcommand{\bm}{\mathbf{m}}
\newcommand{\ba}{\mathbf{a}}
\newcommand{\bb}{\mathbf{b}}

\newcommand{\pr}{\emph{Proof.}\ }
\newcommand{\wm}{w_{\textrm{max}}}
\newcommand{\OO}{\mathcal{O}}

\newcommand{\al}{\alpha}
\newcommand{\bt}{\beta}

\newcommand{\e}{\epsilon}

\newcommand{\om}{\varpi}
\newcommand{\ta}{\tau}

\newcommand{\s}{\sigma}

\newcommand{\F}{k}

\newcommand{\cf}{\emph{cf.}}
\newcommand{\ca}{\textrm{can}}
\newcommand{\oc}{\textrm{op$\,$can}}
\newcommand{\RR}{\Gamma}

\newcommand{\M}{\mathcal{M}}

\newcommand{\ka}{\kappa}

\newcommand{\Fl}{F\ell}

\newcommand{\E}{\mathcal{E}}
\newcommand{\T}{\mathcal{T}}
\newcommand{\g}{\gamma}
\newcommand{\mS}{\mathcal{S}}

\newcommand{\proj}{\operatorname{pr}}
\newcommand{\id}{\operatorname{id}}

\newcommand{\D}{\mathcal{D}}

\newcommand{\ad}[1]{\overset{#1}{\textrm{---}}}

\renewcommand{\phi}{\varphi}

\begin{document}

\title[SMT for desingularized Richardson varieties]{Standard Monomial Theory for desingularized Richardson varieties in the flag variety $GL(n)/B$}
\author{Micha\"el Balan}
\address{Universit\'e de Valenciennes\\
Laboratoire de Math\'ematiques\\
Le Mont Houy -- ISTV2\\
F-59313 Valenciennes Cedex 9\\
France}
\email{michael.balan@laposte.net}

\date{March 7, 2013}

\maketitle

\begin{abstract}
We consider a desingularization $\RR$ of a Richardson variety in the variety of complete flags, obtained as a fibre of a projection from a certain Bott-Samelson variety $Z$. For any projective embedding of $Z$ via a complete linear system, we construct a basis of the homogeneous coordinate ring of $\RR$ inside $Z$, indexed by combinatorial objects which we call \emph{$w_0$-standard tableaux}. 
\end{abstract}

\section*{Introduction}

Standard Monomial Theory (SMT) originated in the work of Hodge \cite{Ho}, who considered it in the case of the Grassmannian $G_{d,n}$ of $d$-subspaces of a (complex) vector space of dimension $n$. The homogeneous coordinate ring $\C[G_{d,n}]$ is the quotient of the polynomial ring in the Pl\"ucker coordinates $p_{i_1\dots i_d}$ by the Pl\"ucker relations, and Hodge provided a combinatorial rule to select, among all monomials in the $p_{i_1\dots i_d}$, a subset that forms a basis of $\C[G_{d,n}]$: these (so-called standard) monomials are parametrized by semi-standard Young tableaux. Moreover, he showed that this basis is compatible with any Schubert variety $X\subset G_{d,n}$, in the sense that those basis elements that remain non-zero when restricted to $X$ can be characterized combinatorially, and still form a basis of $\C[X]$. The aim of SMT is then to generalize Hodge's result to any flag variety $G/P$ ($G$ a connected semi-simple group, $P$ a parabolic subgroup): in a more modern formulation, the problem consists, given a line bundle $L$ on $G/P$, in producing a ``nice'' basis of the space of sections $H^0(X,L)$ ($X\subset G/P$ a Schubert variety), parametrized by some combinatorial objects. SMT was developed by Lakshmibai and Seshadri (see \cite{LS1,LS2}) for groups of classical type, and Littelmann extended it to groups of arbitrary type (including in the Kac-Moody setting), using techniques such as the path model in representation theory \cite{L1,L2} and Lusztig's Frobenius map for quantum groups at roots of unity \cite{L3}. Standard Monomial Theory has numerous applications in the geometry of Schubert varieties: normality, vanishing theorems, ideal theory, singularities, and so on \cite{LLM1}.

Richardson varieties, named after \cite{R}, are intersections of a Schubert variety and an opposite Schubert variety inside a flag variety $G/P$. They previously appeared in \cite[Ch.~XIV, \S4]{HP} and \cite{St}, as well as the corresponding open subvarieties in \cite{De}. They have since played a role in different contexts, such as equivariant K-theory \cite{LL}, positivity in Grothendieck groups \cite{B2}, standard monomial theory \cite{BL}, Poisson geometry \cite{GY}, positroid varieties \cite{KLS}, and their generalizations \cite{KLS2,BC}. In particular, SMT on $G/P$ is known to be compatible with Richardson varieties \cite{LL} (at least for a very ample line bundle on $G/P$).

Like Schubert varieties, Richardson varieties may be singular \cite{KL,K,U,Ba}. Desingularizations of Schubert varieties are well known: they are the Bott-Samelson varieties \cite{BS,D,Ha}, which are also used for example to establish some properties of Schubert polynomials \cite{M2}, or to give criteria for the smoothness of Schubert varieties \cite{G1,CC}. An SMT has been developed for Bott-Samelson varieties of type~A in \cite{LM}, and of arbitrary type in \cite{LLM2} using the path model \cite{L1,L2}.

In the present paper, we shall describe a Standard Monomial Theory for a desingularization of a Richardson variety. To be more precise, we introduce some notations. Let $G$ be a connected reductive algebraic group over an algebraically closed field $\F$ of arbitrary characteristic, $B$ a Borel subgroup, and $T \subset B$ a maximal torus. We then have a system of simple roots $\al_1,\dots,\al_l$ ($l$ the semisimple rank of $G$), and simple reflections $s_1,\dots,s_l$ that generate the Weyl group $W=N_G(T)/T$. Denote by $B^-$ the Borel subgroup of $G$ opposite to $B$. The $T$-fixed points of $G/B$ are of the form $wB$ for $w\in W$. The Richardson variety $X_w^v \subset G/B$ is the intersection of the direct Schubert variety $X_w=\overline{B.wB}$ with the opposite Schubert variety $X^v=\overline{B^-.vB}=w_0X_{w_0v}$. Fix a reduced decomposition $w=s_{i_1}\dots s_{i_d}$ and consider the Bott-Samelson desingularization $Z=Z_{i_1\dots i_d}(eB)\to X_w$, and similarly $Z'=Z_{i_ri_{r-1}\dots i_{d+1}}(w_0B)\to X^v$ for a reduced decomposition $w_0v=s_{i_{r}}s_{i_{r-1}}\dots s_{i_{d+1}}$. Then the fibred product $Z\times_{G/B}Z'$ has been considered as a desingularization of $X_w^v$ in \cite{B3}, but for our purposes, it will be more convenient to realize it as the fibre $\RR_\bi$ ($\bi=i_1\dots i_di_{d+1}\dots i_r$) of the projection $Z_\bi=Z_\bi(eB)\to G/B$ over $w_0B$ (see Section~\ref{section-Desingularized} for the precise connection between those two constructions).

In \cite{LM,LLM2}, Lakshmibai, Littelmann, and Magyar define a family of line bundles $L_{\bi,\bm}$ ($\bm=m_1\dots m_r\in\Z_{\ge0}^r$) on $Z_\bi$ (they are the only globally generated line bundles on $Z_\bi$, as pointed out in \cite{LT}), and give a basis for the space of sections $H^0(Z_{\bi},L_{\bi,\bm})$. In \cite{LM}, the elements $p_T$ of this basis, called standard monomials, are indexed by combinatorial objects $T$ called standard tableaux: the latter's definition involves certain sequences $J_{1,1} \supset \dots\supset J_{1,m_1}\supset\dots \supset J_{r,1}\supset \dots \supset J_{r,m_r}$ of subwords of $\bi$, called liftings of $T$ (see Section~\ref{section-LM} for precise definitions---actually, two equivalent definitions of standard tableaux are given in \cite{LM}, but we will only use the one in terms of liftings). Note also that $L_{\bi,\bm}$ is very ample precisely when $m_j>0$ for all $j$ (see \cite{LT}, Theorem~3.1), in which case $\bm$ is called regular.

The main result of this paper states that in this case, if $G=GL(n,\F)$, then SMT on $Z_\bi$ is compatible with $\RR_\bi$.
\newpage

\begin{The}\label{main-theorem}
Let $G=GL(n,\F)$. Assume that $\bm$ is regular. With the above notation, the standard monomials $p_T$ such that $(p_T)_{|\RR_\bi} \neq 0$ still form a basis of $H^0(\RR_\bi,L_{\bi,\bm})$.\\
Moreover, $(p_T)_{|\RR_\bi} \neq 0$ if and only if $T$ admits a lifting $J_{1,1}\supset \dots\supset J_{r,m_r}$ such that each subword $J_{k,m}$ contains a reduced expression of $w_0$.
\end{The}

We prove this theorem in three steps.
\begin{enumerate}
 \item Call $T$ (or $p_T$) $w_0$-standard if the above condition on $(J_{k,m})$ holds. We prove by induction over $M=\sum_{j=1}^r m_j$ that the $w_0$-standard monomials $p_T$ are linearly independent on $\RR_\bi$. (Here the assumption that $\bm$ is regular is not necessary.)
 \item In the regular case, we prove that a standard monomial $p_T$ does not vanish identically on $\RR_\bi$ if and only if it is $w_0$-standard, using the combinatorics of the Demazure product (see Definition~\ref{Demazure}). It follows that $w_0$-standard monomials form a basis of the homogeneous coordinate ring of $\RR_\bi$ (when $\RR_\bi$ is embedded in a projective space via the very ample line bundle $L_{\bi,\bm}$).
 \item We use cohomological techniques to prove that the restriction map 
\[
H^0(Z_{\bi},L_{\bi,\bm}) \to H^0(\RR_{\bi},L_{\bi,\bm})
\]
is surjective. More explicitly, we define a family $(Y_\bi^u)$ of subvarieties of $Z_\bi$ indexed by $S_n$, with the property that $Y_\bi^e=Z_\bi$ and $Y_\bi^{w_0}=\RR_\bi$. We construct a sequence in $S_n$, $e=u_0<u_1<\dots<u_N=w_0$, such that for every $t$, $Y_\bi^{u_{t+1}}$ is defined in $Y_\bi^{u_t}$ by the vanishing of a single Pl\"ucker coordinate $p_\ka$, in such a way that each restriction map $H^0(Y_\bi^{u_t},L_{\bi,\bm}) \to H^0(Y_\bi^{u_{t+1}},L_{\bi,\bm})$ can be shown to be surjective using vanishing theorems (Corollary~\ref{Hypersurfaces} and Theorem~\ref{cohomology}). This shows that the $w_0$-standard monomials span $H^0(\RR_\bi,L_{\bi,\bm})$. 
\end{enumerate}

Note that alternate bases of $H^0(Z_\bi,L_{\bi,\bm})$ for certain pairs $(\bi,\bm)$ have been constructed in \cite{Ta}, and the fibred products $Z\times_{GL(n,\C)/B}Z'$ have been studied from this point of view in \cite{FK}.~\\

Sections are organized as follows: in Section~\ref{section-Desingularized}, we first fix notation and recall information on Bott-Samelson varieties $Z_{\bi}$, and then show that the fibre $\RR_\bi$ of $Z_\bi\to G/B$ over $w_0B$ is a desingularization of the Richardson variety $X_w^v$; this fact is most certainly known to experts, but has not, to our knowledge, appeared in the literature. In Section~\ref{section-LM}, we recall the main results about SMT for Bott-Samelson varieties from \cite{LM}, in particular the definition of standard tableaux. In Section~\ref{section-linear_independence}, we define $w_0$-standard monomials and we prove that they are linearly independent in $\RR_\bi$. In Section~\ref{section-non_vanishing}, we prove that when $\bm$ is regular, a standard monomial does not vanish identically on $\RR_\bi$ if and only if it is $w_0$-standard. We prove in Section~\ref{section-basis} that $w_0$-standard monomials generate the space of sections $H^0(\RR_\bi,L_{\bi,\bm})$.\\

\textbf{Acknowledgements.} I would like to thank Christian~Ohn for helpful discussions, as well as Michel~Brion, St\'ephane~Gaussent, and the referees for their careful reading and valuable remarks. In particular, Michel~Brion pointed out gaps in the original proofs of Proposition~\ref{prop:local} and Theorem~\ref{cohomology}.

\section{Desingularized Richardson varieties}\label{section-Desingularized}

The notations are as in the Introduction. In addition, if $k,l\in\Z$, then we denote by $[k,l]$ the set $\{k,k+1,\dots,l\}$, and by $[l]$ the set $[1,l]$. 

We first recall a number of results on Bott-Samelson varieties (see \emph{e.g.}\ \cite{M2}). Although we mainly work in type~A in the sequel, the constructions and results in the present section are given for an arbitrary connected reductive algebraic group $G$.

\begin{Def}
Denote by $l$ the semisimple rank of $G$. For $i \in [l]$, we denote by $\widehat{P}_i$ the minimal parabolic subgroup associated to the simple root $\al_i$.
Two elements $g_1B,g_2B$ in $G/B$ are called \emph{$i$-adjacent} if $g_1\widehat{P}_i=g_2\widehat{P}_i$, a situation denoted by $g_1B \ad{i} g_2B$. 
\end{Def}

Consider a word $\bi=i_1\dots i_r$ in $[l]$, with $w(\bi)=s_{i_1}\dots s_{i_r}\in W$ not necessarily reduced. A \emph{gallery of type $\bi$} is a sequence of elements $g_iB$ of the form
\begin{equation}\label{gallery}
g_0B \ad{i_1} g_1B \ad{i_2} \dots \ad{i_r} g_rB.
\end{equation}
For a given $g_0B$, the \emph{Bott-Samelson variety} of type $\bi$ starting at $g_0B$ is the set of all galleries (\ref{gallery}), \emph{i.e.}\ the fibred product
\[
Z_\bi(g_0B)=\{g_0B\} \times_{G/\widehat{P}_{i_1}} G/B \times_{G/\widehat{P}_{i_2}} \dots \times_{G/\widehat{P}_{i_r}} G/B
\]
(a subvariety of $\left(G/B\right)^r$). In particular, $Z_{i_1\dots i_r}(g_0B)$ is a $\PP^1$-fibration over the subvariety $Z_{i_1\dots i_{r-1}}(g_0B)$, which shows by induction over $r$ that Bott-Samelson varieties are smooth and irreducible.

Each subset $J=\{j_1<\dots<j_k\}\subset[r]$ defines a subword $\bi(J)=i_{j_1}\dots i_{j_k}$ of $\bi$. We then write $Z_J(g_0B)$ instead of $Z_{\bi(J)}(g_0B)$, and we view it as the subvariety of $Z_\bi(g_0B)$ consisting of all galleries (\ref{gallery}) such that $g_{j-1}B=g_jB$ whenever $j\not\in J$.

In the sequel, we shall only need galleries starting at $eB$ or at $w_0B$; in particular, we write $Z_\bi=Z_\bi(eB)$.

The (diagonal) $B$-action on $\left(G/B\right)^r$ leaves $Z_\bi$ invariant. In particular, the $T$-fixed points of $Z_\bi$ are the galleries of the form
\[
eB\ad{i_1} u_1B \ad{i_2} u_1u_2B \ad{i_3} \dots \ad{i_r} u_1\dots u_rB,
\]
where each $u_j\in W$ is either $e$ or $s_{i_j}$. This gallery will be denoted $e_J\in Z_\bi$, where $J=\{j\mid u_j=s_{i_j}\}=\{j_1<\dots<j_k\}$.

For $j\in[r]$, we denote by $\proj_j:Z_\bi\to G/B$ the projection sending the gallery (\ref{gallery}) to $g_jB$. Note that $w(\bi(J))=s_{i_{j_1}}\dots s_{i_{j_k}}=u_1\dots u_r$, so $\proj_r(e_J)=u_1\dots u_rB=w(\bi(J))B$. 

When $\bi$ is reduced, \emph{i.e.}\ $w=s_{i_1}\dots s_{i_r}$ is a reduced expression in $W$, $gB$ lies in the Schubert variety $X_w$ if and only if there is a gallery of type $\bi=i_1\dots i_r$ from $eB$ to $gB$, hence the last projection $\proj_r$ takes $Z_{\bi}$ surjectively to $X_w$. Moreover, this surjection is birational: it restricts to an isomorphism over the Schubert cell $C_w=B.wB$. Thus, $\proj_r:Z_\bi\to X_w$ is a desingularization of $X_w$, and likewise for the last projection $Z_\bi(w_0B)\to X^{w_0w}$.

When $\bi$ is not necessarily reduced, $\proj_r(Z_\bi)$ may be described as follows. Recall \cite[Definition-Lemma~1]{LM} that the poset $\{w(\bi(J))\mid J\subset[r]\}$ admits a unique maximal element, denoted by $\wm(\bi)$ (so $\wm(\bi)=w(\bi)$ if and only if $\bi$ is reduced).

\begin{Prop}\label{image}
Let $\bi$ be an arbitrary word. Then $\proj_r(Z_{\bi})$ is the Schubert variety $X_w$, where $w=\wm(\bi)$. 
\end{Prop}

\pr Since $\proj_r(Z_{\bi})$ is $B$-stable, it is a union of Schubert cells. But $Z_{\bi}$ is a projective variety, so the morphism $\proj_r$ is closed, hence $\proj_r(Z_{\bi})$ is a union of Schubert varieties, and therefore a single Schubert variety $X_w$ since $Z_{\bi}$ is irreducible. Moreover, the $T$-fixed points $e_J$ in $Z_\bi$ project to the $T$-fixed points $w(\bi(J))B$ in $X_w$, and all $T$-fixed points of $X_w$ are obtained in this way (indeed, if $vB$ is such a point, then its fibre $\proj_r^{-1}(vB)$ is $T$-stable, so it must contain some $e_J$ by Borel's fixed point theorem). In particular, $wB$ corresponds to a choice of $J\subset\{1,\dots,r\}$ such that $w(\bi(J))$ is maximal, hence the result.~$\square$\\

We now turn to the description of a desingularization of a Richardson variety $X_w^v=X_w\cap X^v$, $v \le w\in W$. Let $Z=Z_{i_1\dots i_d}$ for some reduced decomposition $w=s_{i_1}\dots s_{i_d}$ and $Z'=Z_{i_r\dots i_{d+1}}(w_0B)$ for some reduced decomposition $w_0v=s_{i_r}s_{i_{r-1}}\dots s_{i_{d+1}}$. Since $Z$ desingularizes $X_w$ and $Z'$ desingularizes $X^v$, a natural candidate for a desingularization of $X_w^v$ is the fibred product $Z \times_{G/B} Z'$. However, we wish to see this variety in a slightly different way: an element of $Z\times Z'$ is a pair of galleries
\begin{gather*}
eB \ad{i_1} g_1B \ad{i_2}\dots \ad{i_d} g_dB,
\\
w_0B \ad{i_r} h_{r-1}B \ad{i_{r-1}} \dots \ad{i_{d+1}} h_dB,
\end{gather*}
and it belongs to $Z\times_{G/B}Z'$ when the end points $g_dB$ and $h_dB$ coincide; in this case, by reversing the second gallery, they concatenate to form a longer gallery
\[
eB\ad{i_1} g_1B\ad{i_2} \dots\ad{i_d} g_dB \ad{i_{d+1}} \dots\ad{i_r} w_0B.
\]
Thus, $Z\times_{G/B}Z'$ identifies with the set of all galleries in $Z_\bi=Z_{i_1\dots i_r}$ that end in $w_0B$, \emph{i.e.}\ with the fibre
\[
\RR_\bi=\proj_r^{-1}(w_0B)
\]
of the last projection $\proj_r:Z_\bi\to G/B$. By construction, the $d$th projection $\proj_d$ then maps $\RR_\bi$ onto the Richardson variety $X_w^v$.\\

\begin{Prop}~\label{desingularization}
In the above notation, the $d$th projection $\proj_d:\RR_\bi\to X_w^v$ is a desingularization, \emph{i.e.}\ $\proj_d$ is birational, and the variety $\RR_{\bi}$ is smooth and irreducible.
\end{Prop}
 
I would like to thank the referee for an argument of smoothness and irreducibility that is much simpler than the one in the first version of this paper.\\

\pr  Let $U$ be the maximal unipotent subgroup of $B$. Then $U\to C_{w_0}, u \mapsto uw_0B$ is an isomorphism. Moreover, $\proj_r: Z_\bi\to G/B$ is surjective and $U$-equivariant, so the map
\[
\begin{array}{rcl}
U \times \RR_\bi &\to& \proj_r^{-1}(C_{w_0})\\
(u,\g) & \mapsto & u.\gamma 
\end{array}
\]
is an isomorphism. Hence $Z_\bi$ contains an open set isomorphic to $U\times \RR_\bi$. Since $Z_\bi$ is smooth and irreducible, so is $\RR_\bi$.

Finally, to show that $\RR_\bi\to X_w^v$ is birational, we consider the projections $\proj_d:Z\to X_w$ and $\proj_{r-d}:Z'\to X^v$. Since they are birational, there exist open subsets $U_w\subset X_w$ and $O\subset Z$ isomorphic under $\proj_d$, and open subsets $U^v\subset X^v$ and $O'\subset Z'$ isomorphic under $\proj_{r-d}$. Then the open set $(O\times O') \cap (Z \times_{G/B} Z')$ of $Z \times_{G/B} Z'$ is isomorphic to the open set $U_w \cap U^v$ of $X_w^v$ under $\proj_d:Z\times_{G/B}Z'\to X_w^v$. Since $X_w^v$ and $Z\times_{G/B}Z'\cong\RR_\bi$ are irreducible, these open subsets must be dense. The birationality of $\proj_d:\RR_\bi\to X_w^v$ follows.~$\square$
~\\

\begin{Rem}
In characteristic 0, the smoothness of the fibred product $Z \times_{G/B} Z'$ can alternatively be proved using Kleiman's transversality theorem (\cf\ \cite{H}, Theorem~10.8). This theorem may also be used to prove the irreducibility of $Z \times_{G/B} Z'$, as follows. First, Kleiman's theorem states that every irreducible component of $Z \times_{G/B} Z'$ is of dimension $\dim(Z)+\dim(Z')-\dim(G/B)$. Now consider $\partial Z$ (resp. $\partial Z'$) the union of all Bott-Samelson varieties $X$ with $X \subsetneq Z$ (resp. $X \subsetneq Z'$). Again by Kleiman's theorem, the dimension of $(\partial Z \times_{G/B} Z') \cup (Z \times_{G/B} \partial Z')$ is less than $\dim(Z \times_{G/B}Z')$. So, on one hand, the fibred product $O=(Z \setminus \partial Z) \times_{G/B} (Z' \setminus \partial Z')$ meets each irreducible component of $Z \times_{G/B} Z'$, hence $O$ is dense. On the other hand, $O$ is isomorphic to the open Richardson variety $C_w^v=C_w\cap C^v$, where $C^v$ is the opposite Schubert cell $B^-.vB$, hence $O$ is irreducible, and so is $Z \times_{G/B} Z'$.
\end{Rem}
~\\

For $\bi$ an arbitrary word, we may still consider the variety $\RR_\bi$ of galleries of type $\bi$, beginning at $eB$ and ending at $w_0B$. In general this variety is no longer birational to a Richardson variety. But we still have\\

\begin{Prop}~\label{imageRR}
Let $\bi=i_1\dots i_r$ be an arbitrary word, and consider the projection $\proj_j: \RR_\bi \to G/B$. Then $\proj_j(\RR_\bi)$ is the Richardson variety $X_y^x$ where $y={\wm(i_1\dots i_j)}$ and $x={w_0\wm(i_{j+1}\dots i_r)^{-1}}$. Moreover, $\RR_\bi$ is smooth and irreducible.
\end{Prop}

\pr The variety $\RR_\bi$ is isomorphic to the fibred product 
\[
Z_{i_1\dots i_j}\times_{G/B} Z_{i_r\dots i_{j+1}}(w_0B),
\]
hence 
\[
\begin{split}
\proj_j(\RR_\bi)&=\proj_j(Z_{i_1\dots i_j}) \cap \proj_{r-j}(Z_{i_r\dots i_{j+1}}(w_0B))\\
                &=X_{\wm(i_1\dots i_j)} \cap w_0X_{\wm(i_r\dots i_{j+1})}\\
                &=X_y^x.
\end{split}
\]
Finally, we may prove that $\RR_\bi$ is smooth and irreducible exactly as in the proof of Proposition~\ref{desingularization}.~$\square$
~\\

\begin{Nots}
In subsequent sections, we shall work in type~A. In this setting, we take $G=GL(n,\F)$, which is of semisimple rank $n-1$, $B$ the group of upper triangular matrices, $B^-$ the group of lower triangular matrices, and $T$ the group of diagonal matrices. The variety $G/B$ then identifies with the variety $\Fl(n)$ of complete flags in $\F^n$. The Weyl group $W$ is isomorphic to the group $S_n$ of permutations. For $i\neq j$, we denote by $\al_{i,j}$ the root such that the associated reflection $s_{\al_{i,j}}\in S_n$ is the transposition $(i,j)$. The simple roots are then the $\al_i:=\al_{i,i+1}$. For any element $w\in S_n$, we shall also use the one-line notation $[w(1)\dots w(n)]$.

Let $(e_1,\dots,e_n)$ denote the canonical basis of $\F^n$. The $T$-fixed point $wB$ associated to the permutation $w\in S_n$ is the flag whose $i$th component is the space spanned by $e_{w(1)},\dots,e_{w(i)}$. The flags corresponding to $eB$ and $w_0B$ are respectively denoted by $F_\ca$ and $F_\oc$.
\end{Nots}

\section{Background on SMT for Bott-Samelson varieties}\label{section-LM}

In this section, we recall from \cite{LM} the main definitions and results about Standard Monomial Theory for Bott-Samelson varieties for $G=GL(n,\F)$.

\begin{Defs}
A \emph{tableau} is a sequence $T=t_1\dots t_p$ with $t_j \in [n]$. If $T=t_1\dots t_p$ and $T'=t'_1\dots t'_{p'}$ are two tableaux, then the \emph{concatenation} $T*T'$ is the tableau $t_1\dots t_p\,t'_1\dots t'_{p'}$. We denote by $\emptyset$ the empty tableau, so that $T*\emptyset=\emptyset*T=T$.

A \emph{column $\ka$ of size $i$} is a tableau $\ka=t_1\dots t_i$ with $1\le t_1<\dots<t_i\le n$. The set of all columns of size $i$ is denoted by $I_{i,n}$. The \emph{Bruhat order} on $I_{i,n}$ is defined by
\[
 \ka=t_1\dots t_i\le \ka'=t'_1\dots t'_i \iff t_1\le t'_1,\ \dots, \ t_i\le t'_i.
\]
The symmetric group $S_n$ acts on $I_{i,n}$: if $w\in S_n$ and $\ka=t_1\dots t_i \in I_{i,n}$, then $w\ka$ is the column obtained by rearranging the tableau $w(t_1)\dots w(t_i)$ in an increasing sequence.

For $i \in [n]$, the \emph{fundamental weight column $\om_i$} is the sequence $12\dots i$.
\end{Defs}

We shall be interested in a particular type of tableaux, called standard.

\begin{Defs}
Let $\bi=i_1\dots i_r$, and $\bm=m_1\dots m_r\in \Z_{\ge 0}^r$. A \emph{tableau of shape} $(\bi,\bm)$ is a tableau of the form 
\[
\ka_{1,1}*\dots*\ka_{1,m_1}*\ka_{2,1}*\dots*\ka_{2,m_2}*\dots*\ka_{r,1}*\dots*\ka_{r,m_r} 
\]
where $\ka_{k,m}$ is a column of size $i_k$ for every $k,m$. (If $m_k=0$, there is no column in the corresponding position of $T$.)

A \emph{lifting} of $T$ is a sequence of subwords of $\bi$, 
\[
J_{1,1} \supset \dots \supset J_{1,m_1} \supset J_{2,1} \supset \dots \supset J_{2,m_2} \supset \dots \supset J_{r,1} \supset \dots \supset J_{r,m_r},
\]
such that $J_{k,m} \cap [k]$ is a reduced subword of $\bi$ and $w(\bi(J_{k,m} \cap [k]))\om_{i_k}=\ka_{k,m}$. If such a lifting exists, then the tableau $T$ is said to be \emph{standard}.
\end{Defs}

\begin{Rem}\label{geomlifting}
The last equality in the definition of a lifting may be viewed geometrically as follows. If $J\subset[r]$ and $j\in[r]$, then $\proj_j:Z_\bi\to \Fl(n)$ maps $Z_J\subset Z_\bi$ onto a Schubert variety $X_w\subset \Fl(n)$ (\cf\ Proof of Proposition~\ref{image}). In the notations of Section~\ref{section-Desingularized}, the images of $T$-fixed points of $Z_J$ under $\proj_j$ are of the form $\proj_j(e_K)=e_{u_1\dots u_j}=e_{w(\bi(K\cap[j]))}$ with $K$ running over all subsets of $J$, hence $w=\wm(\bi(J\cap[j]))$. In turn, the image of $\proj_j(Z_J)$ by the projection $\Fl(n) \to G_{i_j,n}$ is equal to the Schubert variety $X_{w\om_{i_j}}$: for $J=J_{k,m}$ in the above lifting, this projection is therefore equal to $X_{\ka_{k,m}}$. We shall follow up on this point of view in Remark~\ref{geomopt}.
\end{Rem}

\begin{Not}
Each column $\ka \in I_{i,n}$ identifies with a weight of $GL(n)$, in such a way that the fundamental weight column $\om_i$ corresponds to the $i$th fundamental weight of $GL(n)$. Therefore, we also denote by $\om_i$ this fundamental weight.
\end{Not}

We recall the Pl\"ucker embedding: given an $i$-subspace $V$ of $\F^n$, choose a basis $v_1,\dots,v_i$ of $V$, and let $M$ be the matrix of the vectors $v_1,\dots,v_i$ written in the basis $(e_1,\dots,e_n)$. We associate to each column $\ka=t_1\dots t_i$ the minor $p_\ka(V)$ of $M$ on rows $t_1,\dots,t_i$. Then the map $p:\ V \mapsto [p_\ka(V)\mid \ka \in I_{i,n}]$ is the Pl\"ucker embedding. 

Let $\pi_i: \Fl(n) \to G_{i,n}$ be the natural projection. We denote by $L_{\om_i}$ the line bundle $(p\circ \pi_i)^*\OO(1)$.

Now consider the tensor product $L_{\om_{i_1}}^{\otimes m_1} \otimes\dots\otimes L_{\om_{i_r}}^{\otimes m_r}$ on $\Fl(n)^r$, and denote by $L_{\bi,\bm}$ its restriction to $Z_\bi\subset\Fl(n)^r$.

\begin{Def}\label{standard_monomial}
To a tableau $T=\ka_{1,1}*\dots*\ka_{1,m_1}*\dots*\ka_{r,1}*\dots*k_{r,m_r}$, one associates the section $p_T=p_{\ka_{1,1}}\otimes \dots \otimes p_{\ka_{1,m_1}} \otimes \dots \otimes p_{\ka_{r,1}} \otimes \dots \otimes p_{\ka_{r,m_r}}$ of $L_{\bi,\bm}$. If $T$ is standard of shape $(\bi,\bm)$, then $p_T$ is called a \emph{standard monomial of shape $(\bi,\bm)$}.
\end{Def}
\vspace*{2mm}

\begin{The}[{\cite{LM}}]~\label{SMT-BS}
\begin{enumerate}
\item The standard monomials of shape $(\bi,\bm)$ form a basis of the space of sections $H^0(Z_{\bi},L_{\bi,\bm})$.
\item For $i>0$, $H^i(Z_{\bi},L_{\bi,\bm})=0$.
\item The variety $Z_\bi$ is projectively normal for any embedding induced by a very ample line bundle $L_{\bi,\bm}$.
\end{enumerate}
\end{The}

\section{Linear independence}~\label{section-linear_independence}

In this section, we define the notion of $w_0$-standard monomials, and then prove that they are linearly independent.

\begin{Defs}~\label{def-w0-standard}
Let $T$ be a standard tableau of shape $(\bi,\bm)$. We say that $T$ (or the monomial $p_T$) is \emph{$w_0$-standard} if there exists a lifting $(J_{k,m})$ of $T$ such that each subword $J_{k,m}$ contains a reduced expression of $w_0$. 

More generally, if $J \subset [r]$ contains a reduced expression for $w_0$, then $\RR_J=Z_J \cap \RR_\bi\neq \emptyset$, and we say that $T$ (or $p_T$) is \emph{$w_0$-standard on $\RR_J$} if there exists a lifting $(J_{k,m})$ of $T$ such that for every $k,m$, $J \supset J_{k,m}$ and $J_{k,m}$ contains a reduced expression of $w_0$. 

Similarly, $T$ (or $p_T$) is said to be \emph{$w_0$-standard on a union $\RR=\RR_{J_1} \cup \dots \cup \RR_{J_k}$} if $T$ is $w_0$-standard on at least one of the components $\RR_{J_1},\dots,\RR_{J_k}$. We then denote by $\mS(\RR)$ the set of all $w_0$-standard tableaux on $\RR$.
\end{Defs}

We need some results about \emph{positroid varieties} (see \cite{KLS}). Let $\pi_i$ be the canonical projection $\Fl(n) \to G_{i,n}$. In general, the projection of a Richardson variety $X_w^v \subset \Fl(n)$ is no longer a Richardson variety. But $\pi_i(X_w^v)$ is still defined inside the Grassmannian $G_{i,n}$ by the vanishing of some Pl\"ucker coordinates. More precisely, consider the set $\M=\{\ka \in I_{i,n} \mid e_{\ka} \in \pi_i(X_w^v)\}$. Then 
\[
\pi_i(X_w^v)=\{V \in G_{i,n} \mid \ka \notin \M \implies p_\ka(V)=0 \}.\qquad(*)
\]
The poset $\M$ is a \emph{positroid} (see the paragraph following Lemma~3.20 in \cite{KLS}), and the variety~$(*)$ is called a \emph{positroid variety}.

\begin{Lem}\label{positroid}
 With the notation above, 
\[
 \M=\{\ka \in I_{i,n}\mid \exists u \in [v,w],\ u\om_i=\ka\}.
\]
\end{Lem}

\pr Let $u \in [v,w]$ and $\ka=u\om_i$. Then $e_u \in X_w^v$, so $e_\ka=\pi_i(e_u)\in\pi_i(X_w^v)$. Hence $\ka \in \M$.

Conversely, let $\ka \in \M$. The fibre $\pi_i^{-1}\{e_\ka\}$ in $X_w^v$ is a non-empty $T$-stable variety, hence, by Borel's fixed point theorem, this variety has a $T$-fixed point $e_u$, $u \in S_n$. It follows that $u \in [v,w]$ and $u\om_i=\ka.~\square$\\

\begin{The}\label{Linear independence}
For every subword $J_1,\dots,J_k$ containing a reduced expression of $w_0$, the $w_0$-standard monomials on the union $\RR=\RR_{J_1}\cup \dots \cup \RR_{J_k}$ are linearly independent.
\end{The}

\pr We imitate the proof of the corresponding proposition for Bott-Samelson varieties appearing in \cite[Section~3.2]{LM}. Let $\T$ be a non-empty subset of $\mS(\RR)$, and assume that we are given a linear relation among  monomials $p_T$ for $T$ in $\T$:
\[
\sum_{T \in \T}a_T\,p_T(\g)=0 \quad \forall \g \in \RR.\qquad(**)
\]
Moreover, we may assume that the coefficients appearing in this relation are all non-zero.
We shall proceed by induction on the length of tableaux, that is, on $M=\sum_{i=1}^r m_i$. 

If $M=1$, then $\bm$ has the form $0\dots 1\dots0$, that is, we have $m_e=1$ for some $e$, and $m_i=0$ for all $i \neq e$. The tableaux $T$ that appear in relation $(**)$ are
of the form $T=\ka_e$, where $\ka_e \in I_{i_e,n}$. If $\g=(F_\ca,F_1,\dots,F_r=F_\oc) \in \RR$ then $p_T(\g)=p_{\ka_e}(F_e)$. Thus, we have a linear relation of Pl\"ucker 
coordinates in a union of Richardson varieties in $\Fl(n)$, hence a linear relation on one of these Richardson varieties. But Standard Monomial Theory for Richardson varieties
(\textit{cf.} \cite{LL}, Theorem~32) shows that such a relation cannot exist.
 
Now assume that $M>1$, and $\bm=0\dots 0\,m_e\dots m_r$ with $m_e>0$. Here, we denote by $\ka_{k,m}^T$ the columns of a tableau $T$. Consider an element $\ka$ minimal among the first columns of the tableaux of $\T$, that is,
\[
 \ka \in \min\{\ka_{e,1}^T\ |\ T \in \T\}.
\]
 
We consider the set $\T(\ka)$ of tableaux $T$ in $\T$ with $\ka_{e,1}^T=\ka$. For every $T \in \T(\ka)$, fix a maximal lifting $J_{e,1}^T \supset
 \dots \supset J_{r,m_r}^T$ containing a reduced expression of $w_0$ and with $J_{e,1}^T$ contained in one of the subwords $J_1,\dots,J_k$, so that $\RR\supset \RR_{J_{e,1}^T}\neq \emptyset$. Thus, we can restrict the relation $(**)$ on 
\[
\RR(\ka)=\bigcup_{T \in \T(\ka)}\RR_{J_{e,1}^T}.
\]

If $T \in \T(\ka)$, then $T=\ka*T'$, and $T'$ is a $w_0$-standard tableau on $\RR(\ka)$ of shape $(\bi,0\dots 0\,m_e-1 \dots m_r)$.

If $T \notin \T(\ka)$, then $\ka_{e,1}^T \nleq \ka$, so $p_{\ka_{e,1}^T}$ vanishes identically on the Schubert variety $X_\ka \subset G_{i_e,n}$, hence on each Schubert variety $X_{\wm(\bi(J_{e,1}^S))}$ for $S \in \T(\ka)$. In particular, $p_{\ka_{e,1}^T}$ vanishes on $\RR(\ka)$, and $p_T$ as well.  

Restrict relation $(**)$ to $\RR(\ka)$: 
\[
p_\ka(\g) \sum_{T\in \T(\ka)} a_{T}\,p_{T'}(\g)=0 \quad \forall \g \in \RR(\ka). 
\]

This product vanishes on each irreducible $\RR_{J_{e,1}^T}$ ($T\in \T(\ka)$). Now, $p_\ka$ does not vanish identically on $\RR(J_{e,1}^T)$. Indeed, we know by Proposition~\ref{imageRR} that $\proj_e(\RR(J_{e,1}^T))$ is the Richardson variety $X_y^x$ with $y=w(\bi(J_{e,1}^T))\ge x$. Since $\ka=y\om_{i_e}$, by Lemma~\ref{positroid}, $p_\ka$ does not vanish identically on $X_y^x$, hence does not vanish identically on $\RR(J_{e,1}^T)$. 

So we may simplify by $p_\ka$ on the irreducible $\RR_{J_{e,1}^T}$, hence a linear relation between $w_0$-standard monomials on $\RR(\ka)$ of shape $(\bi,0\dots 0\,m_e-1\dots m_r)$. By induction over $M$, $a_T=0$ for all $T \in \T(\ka)$: a contradiction.~$\square$

\section{Standard monomials that do not vanish on $\RR_\bi$ are $w_0$-standard}~\label{section-non_vanishing}

In this section, we shall prove that the standard monomials that do not vanish identically on $\RR_\bi$ are $w_0$-standard, provided certain assumptions over $\bm$, which cover the regular case (\emph{i.e.} $m_i>0$ for every $i$). We use a number of combinatorial facts that are true in arbitrary type. We state them in a separate subsection.

\subsection{Combinatorics of the Demazure product}

In this subsection, we assume $W$ to be any finite Weyl group, with Bruhat order $\leq$ and length function $\ell$.

\begin{Lem}[Lifting Property {\cite[Proposition 2.2.7]{BB}}]\label{Z}
Let $s\in W$ be a simple reflection, and $u<w$ in $W$.
\begin{itemize}
 \item If $u<su$ and $w>sw$, then $u \le sw$ and $ su \le w$.
 \item If $u>su$ and $w>sw$, then $su\le sw$.
 \item If $u<su$ and $w<sw$, then $su\le sw$.~$\square$
\end{itemize}
\end{Lem}

We may represent these situations by the pictures below
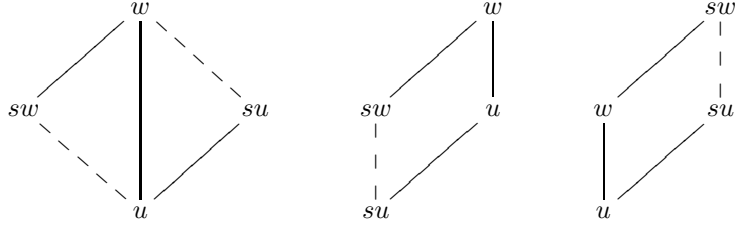
\[
\begin{array}{c@{\hspace{1cm}}c@{\hspace{1cm}}c}
 \xymatrix @R=1cm @C=1cm {
 & w \ar@{-}[ld] \ar@{--}[rd] \ar@{-}[dd]& \\
sw \ar@{--}[rd]& & su \ar@{-}[ld] \\
 & u &  
}&
 \xymatrix @R=1cm @C=1cm {
 & w \ar@{-}[d] \ar@{-}[dl]\\
sw \ar@{--}[d] & u \ar@{-}[dl] \\
su &
}
 & \xymatrix @R=1cm @C=1cm {
 & sw \ar@{-}[dl] \ar@{--}[d]\\
w \ar@{-}[d] & su \ar@{-}[dl] \\
u &
}
\end{array}
\]
 
\begin{Def}\label{Demazure}
Let $x,y \in W$. The Demazure product $x*y$ is the unique maximal element of the poset $\D(x,y)=\{uv\mid u\le x, v\le y\}$.
\end{Def}

\begin{Prop}[{\cite{RS,HLu}}]\label{prop:dblcoset}
Let $x,y \in W$. The double coset $B(x*y)B$ is the unique $B\times B$-double coset that is open in $BxByB$. In particular, $*$ is well-defined and associative.~$\square$ 
\end{Prop}

\begin{Lem}\label{lem:s*}
Let $s$ be a simple reflection, and $x \in W$. Then $x*s=\max(x,xs)$. Similarly, $s*x=\max(x,sx)$.
\end{Lem}

\pr We shall prove that $x*s=\max(x,xs)$, the proof of $s*x=\max(x,sx)$ being similar.
\begin{itemize}
 \item Case~1: $x>xs$. Let $u \le x$. If $us<u$, then $us\le x$. If $us>u$, then by Lemma~\ref{Z}, we have $us\le x$. Hence every element of $\D(x,s)$ is less than or equal to $x$, so $x*s=x=\max(x,xs)$.
 \item Case~2: $x<xs$. Let $u\le x$. If $us<u$, then $us\le xs$. If $us>u$, then by Lemma~\ref{Z}, $us\le xs$. Thus, every element of $\D(x,s)$ is less than or equal to $xs$, so $x*s=xs=\max(x,xs).~\square$ 
\end{itemize}
\vspace*{2mm}

\begin{Lem}~\label{decompose}
 Let $J$ be a subword of $\bi$. For every $k \in [r]$, 
\[
\wm\bigl(\bi(J)\bigr)=\wm\bigl(\bi(J\cap[k])\bigr)*\wm\bigl(\bi(J\cap[k+1,r])\bigr). 
\]
\end{Lem}

\pr Let 
\begin{align*}
w=&\wm\bigl(\bi(J)\bigr)\\ 
x=&\wm\bigl(\bi(J\cap[k])\bigr)\\
y=&\wm\bigl(\bi(J\cap[k+1,r])\bigr)
\end{align*}

Each element $uv$ of $\D(x,y)$ has a decomposition of the form $w\bigl(\bi(K_1)\bigr)w\bigl(\bi(K_2)\bigr)$ with $K_1\subset J\cap[k]$ and $K_2\subset J\cap[k+1,r]$. Hence, 
\[
uv=w\bigl(\bi(K_1\cup K_2)\bigr)\le w,
\] 
so $x*y\le w$.

Conversely, let $K' \subset J$ be such that $w\bigl(\bi(K')\bigr)=w$ is a reduced decomposition. Since 
\[
w=w\bigl(\bi(K'\cap[k])\bigr)\,w\bigl(\bi(K'\cap[k+1,r])\bigr),
\] 
we have $w \in \D(x,y)$, hence $w\le x*y$.~$\square$\\

\begin{Lem}[{\cite[2.2.(4)]{HL}}]~\label{*order}
 If $x'\le x$ and $y'\le y$, then $x'*y'\le x*y$.
\end{Lem}

\pr By Proposition~\ref{prop:dblcoset}, we have
\[
 \begin{split}
  B(x'*y')B \subset \overline{BxB}\ \overline{ByB} \subset \overline{BxByB}=\overline{B(x*y)B},
 \end{split}
\]
hence $x'*y' \le x*y$.~$\square$
\vspace*{2mm}

We shall also need a result due to V.~Deodhar:

\begin{Nots}

Denote by $P$ a parabolic subgroup of $G$ containing $B$, $W_P$ the parabolic subgoup of $W$ associated to $P$, and $W^P$ the set of minimal representatives of the quotient $W/W_P$. This set indexes Schubert varieties and $T$-fixed points of $G/P$. 

Let $\ka \in W^P$ and $w \in W$. We set 
\[
\E(w,\ka)=\{v \in W\mid v\le w,\ v \equiv \ka \bmod W_P\}.
\]
\end{Nots}

\begin{Lem}[{\cite[Lemma~11]{LLM2}}]\label{uniquemax}
Let $\ka\in W^P$, and $w\in W$. If $\E(w,\ka)\neq \emptyset$, then it admits a unique maximal element.~$\square$
\end{Lem}

\begin{Rem}\label{rem:geomuniquemax}
The above lemma admits the following geometric interpretation. Let $q$ be the restriction to $X_w$ of the canonical projection $G/B \to G/P$. Since $q$ is $B$-equivariant, $q^{-1}(X_\ka)$ is a union of Schubert varieties, namely
\[
 q^{-1}(X_\ka)=\bigcup_{v \in \E(w,\ka)} X_v.
\]
By Lemma~\ref{uniquemax}, we conclude that $q^{-1}(X_\ka)$ is a single Schubert variety.\\

There exists a dual version of the lemma: if $w\ge \ka$, then the set
\[
 \{v \in W\mid v\ge w,\ v\equiv \ka \bmod W_P\}
\]
admits a unique minimal element. Equivalently, $q^{-1}(X^\ka)$ is a single opposite Schubert variety.  
\end{Rem}

\begin{Lem}~\label{algorithm}
 Let $w \in W$ and $\ka \in W^P$ be such that $\E(w,\ka)\neq \emptyset$. Consider a simple reflection $s$ such that $sw<w$.
\begin{enumerate}
 \item If $s\ka >\ka$, then $\max\E(w,\ka)=\max\E(sw,\ka)$.
 \item If $s\ka \le \ka$, then $\max\E(w,\ka)=s*\max\E(sw,s\ka)$.
\end{enumerate}
\end{Lem}

\pr Let $u=\max\E(w,\ka)$.
\begin{itemize}
\item Case~1: assume that $su>u$. Then, by Lemma~\ref{Z},
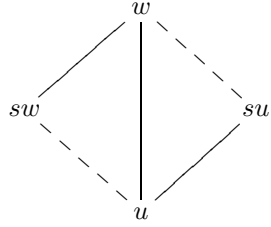
\[
 \xymatrix @R=1cm @C=1cm {
 & w \ar@{-}[dl] \ar@{-}[dd] \ar@{--}[dr] & \\
sw \ar@{--}[dr] & & su \ar@{-}[dl]\\
 & u & 
}
\]
we have $u\le sw$ and $su\le w$. Hence $s\ka\ge \ka$, but by maximality of $u$, $su \notin \E(w,\ka)$, hence $s\ka>\ka$. Since $u\le sw$,  $u \in \E(sw,\ka)$, so 
\[
u\le \max\E(sw,\ka) \le \max\E(w,\ka)=u. 
\]
This proves the part $(1)$ of the lemma.\\
\item Case~2: $su<u$. Then $s\ka\le \ka$, and by Lemma~\ref{Z},
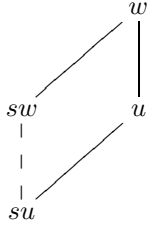
\[
 \xymatrix @R=1cm @C=1cm {
 & w \ar@{-}[dl] \ar@{-}[d]\\
sw \ar@{--}[d] & u \ar@{-}[dl]\\
su & 
}
\]
we have $su\le sw$, so $su \in \E(sw,s\ka)$, hence $su\le v=\max\E(sw,s\ka)$. We distinguish two subcases:
\begin{itemize}
\item Subcase~1: $s\ka<\ka$. Then $sv>v$. Since we also have $su<u$, it follows from Lemma~\ref{Z} that $v\le su$. Similarly, $sv>v$, together with $sw<w$ imply that $sv\le w$, so $sv \in \E(w,\ka)$, hence $sv \le u$. By Lemma~\ref{Z}, we have $v\le su$. So $v=su$, or equivalently 
\[
u=sv=\max(v,sv)=s*v.
\]
\item Subcase~2: $s\ka=\ka$. 
\begin{itemize}
\item If $u\le sw$, then $u \in \E(sw,\ka)$, so $u \le v$. But $v \le u$, so $u=v$. 
\item If $u \nleq sw$, then $su \in \E(sw,\ka)$, so $su\le v\le u$. In other words, $v \in \{u,su\}$. 
\end{itemize}
In each of these two situations, we have $u=v$ or $u=sv$. But, if $sv>v$ then $u \neq v$ (since $su<u$), so $u=sv=\max(v,sv)=s*v$. If $sv<v$, then $u \neq sv$, so 
\[
u=v=s*v.~\square
\]
\end{itemize}
\end{itemize}

Let $w=s_{i_1}\dots s_{i_j}$ be a reduced expression. The lemma above gives an algorithm to find a reduced expression of $u=\max\E(w,\ka)$, say $u=w\bigl(\bi(J)\bigr)$, with $J\subset[j]$:  let $s=s_{i_1}$, and compare $s\ka$ with $\ka$. 
\begin{itemize}
\item If $s\ka>\ka$, then $u=\max\E(sw,\ka)$. 
\item If $s\ka\le \ka$, then $u=s*\max\E(sw,s\ka)$. 
\end{itemize}
We then compute $\max\E(sw,s\ka)$ or $\max\E(sw,\ka)$ in the same way, using the decomposition $sw=s_{i_2}\dots s_{i_j}$.\\

\begin{Ex}
If $W=S_4$ and $P=S_2\times S_2$, then $W^P$ identifies with the set $\{ij\mid1\le i<j\le4\}$. Take $w=[4231]=s_1s_2s_3s_2s_1$ and $\ka=13$. We shall compute $u=\max\E(w,\ka)$ with the previous algorithm. Note that $\ka\le 24=w\om_2$, hence $\E(w,\ka)\neq \emptyset$.
\begin{itemize}
 \item $s_1\ka=23>\ka$, so $u=\max\E(s_2s_3s_2s_1,\ka)$,
 \item $s_2\ka=12\le k$, so $u=s_2*\max\E(s_3s_2s_1,12)$,
 \item $s_3(12)=12$, so $u=s_2*s_3*\max\E(s_2s_1,12)$,
 \item $s_2(12)=13>12$, so $u=s_2*s_3*\max\E(s_1,12)$,
 \item $s_1(12)=12$, so $u=s_2*s_3*s_1*\max\E(e,12)$.
\end{itemize}
Now, $\max(e,12)=e$, so $u=s_2*s_3*s_1=s_2s_3s_1=[3142]$.
\end{Ex}

\begin{Lem}[{\cite[Proposition~2.4.4]{BB}}]\label{factorization}
Let $\ka \in W^P$. The set $\{v\in W\mid v\equiv \ka \bmod W_P\}$ admits a unique minimal element $u$. Moreover, if $v\in W$ satisfies $v\equiv \ka \bmod W_P$, then $v$ admits a unique factorization $v=uv'$ with $v'\in W_P$. This factorization is length-additive, in the sense that $\ell(v)=\ell(u)+\ell(v')$.~$\square$
\end{Lem}

\begin{Not}
We denote by $P_i$ the maximal parabolic subgroup of $G$ containing $B$ associated to the $i$th simple root $\al_i$.
\end{Not}

\begin{Lem}\label{lem_proof}
Denote by $u_d$ the minimal element of $W$ such that $u_d\equiv w_0 \bmod W_{P_d}$. Let $w\ge u$, and $\ka\in W^{P_i}$ for an arbitrary $i$ such that $\E(w,\ka)\neq \emptyset$. Assume that $x=\max\E(w,\ka)\ge u$. Then
\[
 \forall v\ge u,\ \E(v,\ka)\neq \emptyset \implies \max \E(v,\ka)\ge u.
\]
\end{Lem}

\pr Since $v\ge u$, we have $v\equiv w_0\bmod W_{P_d}$, hence by Lemma~\ref{factorization}, $v=uv'$ with $v'$ in $W_{P_d}$. Moreover, this decomposition is length-additive, so if $u=s_{i_1}\dots s_{i_j}$ and $v'=s_{i_{j+1}}\dots s_{i_k}$ are reduced expressions, then $s_{i_1}\dots s_{i_j} s_{i_{j+1}}\dots s_{i_k}$ is a reduced expression of $
v$. Similarly, we decompose $x=ux'$ with $x'\in W_{P_d}$. We then obtain 
\[
 x>s_{i_1}x>\dots>s_{i_j}\dots s_{i_1}x=x',
\]
hence
\[
 \ka\ge s_{i_1}\ka\ge \dots \ge s_{i_j}\dots s_{i_1}\ka.
\]
Now, we apply the procedure described after Lemma~\ref{algorithm} for the decomposition $v=s_{i_1}\dots s_{i_j} s_{i_{j+1}}\dots s_{i_k}$. The above inequalities show that $\max\E(v,\ka)$ is of the form $s_{i_1}*\dots*s_{i_j}*z$. But, by Lemma~\ref{*order}, we have 
\[
\begin{split}
s_{i_1}*\dots*s_{i_j}*z &\ge s_{i_1}*\dots *s_{i_j}\\
                        &\ge s_{i_1}\dots s_{i_j}\\
                        &\ge u.~\square
\end{split}
\]

\subsection{Optimal liftings of a standard tableau}

We now assume again that our Weyl group $W$ is of type~A, \emph{i.e.}\ $W=S_n$.

Let $T$ be a standard tableau of shape $(\bi,\bm)$, and $e$ be the least integer such that $m_e\neq 0$, so $\bm=0\dots0\,m_e\dots m_r$.
We give the construction of a particular type of liftings of $T$ (called optimal), in light of the following

\begin{Rem}\label{geomopt}
Let $(K_{k,m})$ be an arbitrary lifting of $T$ and set 
\[
w_{k,m}=w\bigl(\bi(K_{k,m} \cap [k])\bigr), 
\]
so that $w_{k,m}\om_k=\ka_{k,m}$. By Remark~\ref{geomlifting}, $\proj_k(Z_{K_{k,m}})=X_{w_{k,m}}$, with the following consequences.
\begin{itemize}
 \item For each $k$, $K_{k,1}\supset\dots\supset K_{k,m_k}$ yields $w_{k,1}\ge\dots\ge w_{k,m_k}$.
 \item Let $l$ be the least integer such that $l>k$ and $m_l\neq 0$. Then $K_{k,m_k} \supset K_{l,1}$ yields $\proj_l(Z_{K_{l,1}}) \subset \proj_l(Z_{K_{k,m_k}})$, hence
\[
 w\bigl(\bi(K_{l,1}\cap[l])\bigr)\le\wm\bigl(\bi(K_{k,m_k}\cap[l])\bigr).
\]
By Lemma~\ref{decompose}, 
\[
\wm\bigl(\bi(K_{k,m_k}\cap[l])\bigr)=w\bigl(\bi(K_{k,m_k}\cap[k])\bigr)*\wm\bigl(\bi(K_{k,m_k}\cap [k+1,l])\bigr).
\] 
So
\[
w_{l,1}\le w_{k,m_k}*\wm\bigl(\bi(K_{k,m_k}\cap[k+1,l])\bigr). 
\]
\end{itemize}
\end{Rem}~\\

We now construct elements $v_{k,m}\in S_n$ inductively, as follows. At the first step, consider the set
\[
\E(\wm(i_1\dots i_e),\ka_{e,1}).
\]
Since it contains $w_{e,1}$, it is nonempty, so it has a maximal element $v_{e,1}$, which is unique thanks to Lemma~\ref{uniquemax}. Now assume that $v_{k,m}\ge w_{k,m}$ has already been constructed. We then proceed in the same way to construct $v_{k,m+1}$ (if $m<m_k$) or $v_{l,1}$ (if $m=m_k$, and $l>k$ is the least integer such that $m_l\neq 0$):
\begin{itemize}
\item If $m<m_k$, then the set $\E(v_{k,m},\ka_{k,m+1})$ is nonempty (since it contains $w_{k,m+1}$), so let $v_{k,m+1}$ be its unique maximal element.
\item If $m=m_k$, then let $v'_{k,m}=v_{k,m}*\wm\bigl(i_{k+1}\dots i_l\bigr)$. By Lemma~\ref{*order}, 
\[
w_{k,m}*\wm\bigl(\bi(K_{k,m}\cap[k+1,l])\bigr)\le v'_{k,m}.
\]
Thus, by Remark~\ref{geomopt} the set $\E(v'_{k,m},\ka_{l,1})$  contains $w_{l,1}$, so it is non-empty. Let $v_{l,1}$ be its unique maximal element.
\end{itemize}
\vspace*{2mm}

\begin{Rem}
Although the \emph{existence} of the $v_{k,m}$ depends on that of the $w_{k,m}$ (\emph{i.e.}\ on the existence of a lifting of the tableau $T$), the \emph{values} of the $v_{k,m}$ only depend on the tableau $T$ itself.
\end{Rem}

Next, we construct subsets $E_{k,m}\subset[k]$, again inductively. Since 
\[
v_{e,1}\le \wm(i_1\dots i_e),
\]
choose $E_{e,1}\subset\{i_1\dots i_e\}$ such that $v_{e,1}$ admits a reduced expression of the form $\bi(E_{e,1})$. If $E_{k,m}$ such that $v_{k,m}=w\bigl(\bi(E_{k,m})\bigr)$ has already been constructed, then define $E_{k,m+1}$ (if $m<m_k$) or $E_{l,1}$ (if $m=m_k$) as follows:
\begin{itemize}
\item If $m<m_k$, then $v_{k,m+1}\le v_{k,m}=w\bigl(\bi(E_{k,m})\bigr)$, so choose $E_{k,m+1}\subset E_{k,m}$ such that $v_{k,m+1}$ admits a reduced expression of the form $\bi(E_{k,m+1})$.
\item If $m=m_k$, then by Lemma~\ref{decompose} 
\[
v_{l,1}\le v'_{k,m_k}=\wm\bigl(\bi(E_{k,m_k}\cup\{k+1,\dots,l\})\bigr),
\]
so choose $E_{l,1}\subset E_{k,m_k}\cup\{k+1,\dots,l\}$ such that $v_{l,1}$ admits a reduced expression of the form $\bi(E_{l,1})$.
\end{itemize}
\vspace{2mm}

\begin{Def}~\label{maxlift}
With the above notation, set $J_{k,m}=E_{k,m}\cup[k+1,r]$ for each $k,m$. We will call $(J_{k,m})$ an \emph{optimal lifting} of the tableau $T$.
\end{Def}

\begin{Rem}
An optimal lifting is not unique. However, while it depends on the choice of reduced expressions for the $v_{k,m}$, it is still independent on the choice of the initial lifting $(K_{k,m})$.
\end{Rem}
\vspace*{2mm}

\begin{Not}
 For $k \in [n-1]$, let $j_k$ be the greatest integer such that $i_{j_k}=k$.
\end{Not}

\begin{The}\label{non_vanishing}
Assume that for every $k$, $m_{j_k}>0$. Then the standard monomials $p_T$ of shape $(\bi,\bm)$ that do not vanish identically on $\RR_\bi$ are $w_0$-standard.
\end{The}

\pr Consider an optimal lifting $(J_{k,m})$ of $T$. Let $(F_\ca,F_1,\dots,F_r) \in \RR_\bi$ be a gallery such that $p_T(F_\ca,F_1,\dots,F_r) \neq 0$. By definition of $j_k$, the flags $F_{j_k}$ and $F_\oc$ share the same $k$-subspace, which then is the $T$-fixed point $\langle e_n,\dots,e_{n-k+1}\rangle$. Hence, $\ka_{j_k,1}=\dots=\ka_{j_k,m_{j_k}}=w_0\om_k$.

Arrange the integers $j_1,\dots,j_{n-1}$ in an increasing sequence: $j_{l_1}<\dots<j_{l_{n-1}}$.
 
We shall prove that if $k>j_l$, then $v_{k,m}\ge u_l$. Since $p_T(F_\ca,F_1,\dots,F_r)\neq 0$, we have $p_{\ka_{k,m}}(F_k)\neq 0$, hence $p_{\ka_{k,m}}$ does not vanish identically on the Richardson variety $X_w^v$, where $w=\wm(i_1\dots i_k)$ and $v=w_0(\wm(i_{k+1}\dots i_r))^{-1}$. This means that $p_{\ka_{k,m}}$ does not vanish identically on the positroid variety $\pi(X_w^v)$, where $\pi: \Fl(n) \to G_{i_k,n}$. By Lemma~\ref{positroid}, there exists $u \in [v,w]$ such that $u\om_{i_k}=\ka_{k,m}$. It follows that the maximal element $x_l$ of $\E(w,\ka_{k,m})$ is greater than $u$. But, since $k>j_l$, a reduced expression of $w_0v^{-1}$ consists of letters taken from $i_{k+1}\dots i_r$, where no $l$ appears. Thus, $w_0v^{-1}\om_{l}=\om_{l}$, so $v\om_l=w_0\om_{l}$, that is, $v\ge u_l$. Hence $x_1 \ge u \ge v\ge u_l$. We then conclude with Lemma~\ref{lem_proof}.

Now, we consider subwords $J_{k,m}$ with $k\le j_{l_1}$. In this case, $k\le j_t$ ($t\ge 1$), so we have the inequalities $w\bigl(\bi(J_{j_t,1} \cap [j_t])\bigr)\le \wm\bigl(\bi(J_{k,m})\bigr)$, hence 
\[                                                                                                        
\wm\bigl(\bi(J_{k,m})\bigr)\om_t\ge w\bigl(\bi(J_{j_t,1} \cap [j_t])\bigr)\om_t=\ka_{j_t,1}=w_0\om_t,                                                                                                                                                                                  \]
\emph{i.e.} $\wm\bigl(\bi(J_{k,m})\bigr)\om_t=w_0\om_t$. So $\wm\bigl(\bi(J_{k,m})\bigr)=w_0$.

If $j_{l_t}<k\le j_{l_{t+1}}$, then we have, in one hand, $\wm\bigl(\bi(J_{k,m})\bigr)\ge w\bigl(\bi(J_{j_{l_p},1})\bigr)$ for every $p\ge t+1$, so $\wm\bigl(\bi(J_{k,m})\bigr)\om_{l_p}\ge \ka_{l_p,1}=w_0\om_{l_p}$, hence $\wm\bigl(\bi(J_{k,m})\bigr)\om_{l_p}=w_0\om_{l_p}$. On the other hand, $\wm\bigl(\bi(J_{k,m})\bigr)\ge v_{k,m} \ge u_{l_q}$ for every $q \le t$, hence $\wm\bigl(\bi(J_{k,m})\bigr)\om_{l_q}=w_0\om_{l_q}$. It follows that $\wm\bigl(\bi(J_{k,m})\bigr)=w_0$.~$\square$

\begin{Rem}\label{counter_example}
The assumption  $m_{j_k}>0$ for every $k$ is necessary: consider the standard tableau $T=123*13*3*134*24*124$ of shape $(3213233213,1111110000)$. One may check that $p_T$ does not vanish identically on $\RR_\bi$. However, an optimal lifting of $T$ is given by
\[
 \begin{array}{r @{=\{} c@{}c@{}c@{}c@{}c@{}c@{}c@{}c@{}c@{\}}}
      J_{1,1}    &2,&3,&4,&5,&6,&7,&8,&9,&10 \\[1mm]
      J_{2,1}    &2,&3,&4,&5,&6,&7,&8,&9,&10 \\[1mm]
      J_{3,1}    &2,&3,&4,&5,&6,&7,&8,&9,&10 \\[1mm]
      J_{4,1}    &2,&3,&4,&5,&6,&7,&8,&9,&10 \\[1mm]
      J_{5,1}    &  &3,&4,&5,&6,&7,&8,&9,&10 \\[1mm]
      J_{6,1}    &  &3,& &  &6, &7,&8,&9,&10 \\[1mm]
\end{array}
\]
and we have $\wm(\bi(J_{6,1}))=s_1s_3s_2s_1s_3=[4231]\neq w_0$, hence $T$ is not $w_0$-standard.
\end{Rem}

As an immediate consequence of Theorems~\ref{Linear independence} and \ref{non_vanishing}, we have 

\begin{Cor}\label{cor-basis_coord_ring}
Assume that $\bm$ is regular, so we embed $\RR_\bi$ in the projective space $\PP(H^0(\RR_\bi,L_{\bi,\bm})^*)$. A basis of the homogeneous coordinate ring of $\RR_\bi$ is given by the $w_0$-standard monomials of shape $(\bi,d\bm)$ for all $d \in \N$.~$\square$
\end{Cor}

\section{Basis of $H^0(\RR_{\bi},L_{\bi,\bm})$}\label{section-basis}

Assume that $\bm$ is regular. We shall prove that the $w_0$-standard monomials of shape $(\bi,\bm)$ form a basis of the space of sections $H^0(\RR_\bi,L_{\bi,\bm})$. By Theorems~\ref{Linear independence} and \ref{non_vanishing}, we just have to show that the restriction map 
\[
H^0(Z_\bi,L_{\bi,\bm}) \to H^0(\RR_\bi,L_{\bi,\bm})
\]
is surjective.
The idea is to find a sequence of varieties $\bigl(Y_\bi^{u_t}\bigr)$, parametrized by $u_t \in S_n$, such that
\begin{itemize}
 \item $Y_\bi^{u_0}=Z_\bi$ and $Y_\bi^{u_N}=\RR_\bi$,
 \item $Y_\bi^{u_{t+1}}$ is a hypersurface in $Y_\bi^{u_t}$,
 \item each restriction map $H^0(Y_\bi^{u_{t+1}},L_{\bi,\bm}) \to H^0(Y_\bi^{u_t},L_{\bi,\bm})$ is surjective. 
\end{itemize}
  
\begin{Ex}\label{hypersurface_ex}

Let $n=4$ and $\bi=123212312$. Consider the following reduced expression
\[
 w_0=s_1s_2s_1s_3s_2s_1=s_{a_6}s_{a_5}\dots s_{a_1},
\]
and set
\[
 \left\{\begin{array}{lc}
         u_0=e, & \\
         u_{t+1}=s_{a_{t+1}}u_t & \forall t\ge 0.
        \end{array}\right.
\]
The sequence $(u_{t})$ is increasing, and $u_6=w_0$. Thus, we obtain a sequence of opposite Schubert varieties
\[
 \Fl(n)=X^{u_0} \supset X^{u_1} \supset \cdots \supset X^{u_6}=\{F_\oc\}.
\]

Let $F=(F^1\subset F^2 \subset F^3 \subset F^4=\F^4)$ be a flag.
\begin{itemize}
\item We have the equivalence 
\[
\begin{split}
F\in X^{s_1} &\iff F^1 \in \langle e_2,e_3,e_4\rangle\\
             &\iff p_1(F^1)=0.
\end{split}
\]
\item Assume $F \in X^{u_1}$. Then 
\[
\begin{split} 
F \in X^{s_2s_1} & \iff F^1 \in \langle e_3,e_4\rangle\\
                 & \iff p_2(F^1)=0 \textrm{,\ since we already know that $p_1(F)=0$}.
\end{split}
\]
\item Similarly, $F\in X^{u_3}\iff F\in X^{u_2}\text{\ and\ }p_3(F^1)=0$.
\item $F \in X^{u_4}\iff F\in X^{u_3}\text{\ and\ }p_{14}(F^2)=0$.
\item $F\in X^{u_5}\iff F\in X^{u_4}\text{\ and\ }p_{24}(F^2)=0$.
\item $F\in X^{u_6}\iff F\in X^{u_5}\text{\ and\ }p_{134}(F^3)=0$.
\end{itemize}

We then set $Y_\bi^{u_t}=\proj_9^{-1}(X^{u_t})$. Thus, $Y_\bi^{u_0}=Z_\bi$, $Y_\bi^{u_6}=\RR_\bi$. Moreover, 
\[
\g \in Y_\bi^{u_{t+1}} \iff \g \in Y_\bi^{u_t} \text{\ and\ } p_{\ka_t}(\g)=0,
\]
where we view $\ka_t$ as a tableau of shape $(123212312,\ba'_t)$, with
\[
 \ba'_1=\ba'_2=\ba'_3=000000010,\ \ba'_4=\ba'_5=000000001 ,\ \ba'_6=000000100.
\]
\end{Ex}
\vspace*{2mm}

This example leads us to work with the following varieties. Consider the last projection $\proj_r: Z_{\bi} \to \Fl(n)$. Fix $u\in S_n$ and a reduced decomposition 
\[
w_0u=s_{k_l}s_{k_{l-1}}\dots s_{k_1}.
\]                                                                                                                                                                                                                                                                                                                                      
Consider the opposite Schubert variety $X^u\subset \Fl(n)$ and set
\[
Y_\bi^u=\proj_r^{-1}(X^u)\subset Z_\bi.
\]
In particular, $Y_\bi^e=Z_\bi$ and $Y_\bi^{w_0}=\RR_\bi$.

\begin{Prop}\label{w0-vs-s}
The variety $Y_\bi^u$ is irreducible, and if $\bi'=i_1\dots i_rk_1\dots k_l$, then the projection $\Fl(n)^{r+l}\to \Fl(n)^r$ onto the $r$ first factors restricts to a morphism
\[
\phi: \RR_{\bi'} \to Y_{\bi}^u
\]
that is birational, hence surjective.
\end{Prop}

\pr Recall that a flag $F$ lies in $X^u$ if and only if it can be connected to $F_\oc$ by a gallery of type $k_1\dots k_l$. Hence $Y_\bi^u$ consists of all galleries
\[
F_\ca \ad{i_1} F_1\ad{i_2} \dots\ad{i_r} F_r
\]
that can be extended to a gallery of the form
\[
F_\ca\ad{i_1} F_1\ad{i_2}\dots\ad{i_r} F_r\ad{k_1}\dots\ad{k_l} F_\oc.
\]
Thus, $\phi$ indeed takes values in $Y_\bi^u$ and is surjective. The irreducibility of $Y_\bi^u$ follows.

Moreover, in the diagram
\[
\xymatrix{
\qquad\qquad\RR_{\bi'}\cong Z_\bi\times_{\Fl(n)}Z_{k_l\dots k_1}(F_\oc)\ar[d]\ar[r]&Z_{k_l\dots k_1}(F_\oc)\ar[d]^{\proj_l}
\\
Z_\bi\ar[r]_{\proj_r}&\Fl(n),
}
\]
$\proj_l$ is an isomorphism over $C^u$, and the morphism $\id\times(\proj_l^{-1}\circ\proj_r)$ from $\proj_r^{-1}(C^u)\subset Y_\bi^u$ to $Z_\bi\times Z_{k_l\dots k_1}(F_\oc)$ is an inverse of $\phi$ over $\proj_r^{-1}(C^u)$, hence $\phi$ is birational.~$\square$\\

\begin{Cor}~\label{imageYiu}
Take the notations of the previous proposition, and consider the $j$th projection $\proj_j: Z_\bi \to \Fl(n)$. Then $\proj_j(Y_\bi^u)$ is the Richardson variety $X_y^x$ for $y={\wm(i_1\dots i_j)}$ and $x={w_0\wm(i_{j+1}\dots i_rk_1\dots k_l)^{-1}}$. 
\end{Cor}

\pr Note that $\proj_j(Y_\bi^u)=\proj_j(\phi(\RR_{\bi'}))=\proj_j(\RR_{\bi'})$ since $j \in [r]$. Proposition~\ref{imageRR} leads to the result.~$\square$
~\\

\begin{Nots}
As in Example~\ref{hypersurface_ex}, consider the reduced decomposition 
\[
w_0=s_1(s_2s_1)\dots(s_{n-1}\dots s_1)=s_{a_N}\dots s_{a_1},
\]
and set $u_t=s_{a_t}\dots s_1$, $u_0=e$. 

Consider the sequence of columns $\ka_t$ defined in the following way.
\begin{itemize}
\item The $n-1$ first columns are $\ka_0=1,\ \ka_1=2,\dots,\ \ka_{n-2}=n-1$. 
\item The $n-2$ next columns are $1*n,\ 2*n,\dots,\ n-2*n$. 
\item The $n-3$ next ones are of size 3: $1*w_0\om_2,\ 2*w_0\om_2,\dots,\ n-3*w_0\om_2$.
\item We proceed in the same way for the other columns until we get $\ka_{N-1}=1*w_0\om_{n-2}$.
\end{itemize}

We denote by $b_t$ the size of $\ka_t$, so that $\ka_t=u_t\om_{b_t}$. We set $\ka'_t=u_{t+1}\om_{b_t}$.
\end{Nots}

\begin{Lem}\label{hypersurface_f}
For every $t \in [0,N-1]$, 
\[
X^{u_{t+1}}=\{F \in \Fl(n) \mid F \in X^{u_t} \text{\rm\ and\ }p_{\ka_t}(F)=0\}.
\]

\end{Lem}

\pr We begin by proving the following

{\bf Claim} For every $t$, 

\[
X^{u_{t+1}\om_{b_t}}=\{V \in G_{b_t,n} \mid V \in X^{u_t\om_{b_t}} \text{\ and\ }p_{\ka_t}(V)=0\}.
\]

Indeed, recall that a $b_t$-space $V$ belongs to the opposite Schubert variety $X^{\ka_t}$ if and only if for every $\ka \ngeq \ka_t$, $p_\ka(V)=0$, and similarly for $X^{\ka'_t}$. Thus, we have to describe the set 
\[
 E_t=\{\ka \ngeq \ka'_t\mid \ka \ge \ka_t\}.
\]
We distinguish two cases.
\begin{itemize}
 \item Case~1: $b_{t+1}=b_{t}$. Then $\ka'_t=u_{t+1}\om_{b_t}=\ka_{t+1}$. But $\ka_t$ is of the form $p*w_0\om_{b_t-1}$, and $\ka_{t+1}=(p+1)*w_0\om_{b_t-1}$. So $\ka_t<\ka'_t$, hence $\ka_t \in E_t$. Let $\ka \in E_t$ with $\ka\neq \ka_t$. Then $\ka>\ka_t$, so $\ka \ge \ka_{t+1}$: a contradiction. Hence, the claim is proved in this case.
 \item Case~2: $b_{t+1}=b_t+1$. Then $\ka'_t=u_{t+1}\om_{b_t}=w_0\om_{b_t}=(n-b_t+1)*w_0\om_{b_t-1}$, and $\ka_t=(n-b_t)*w_0\om_{b_t-1}$. Again, $\ka_t \in E_t$. If $\ka\in E_t$ and $\ka>\ka_t$, then $\ka=w_0\om_{b_t}$: a contradiction. This proves the claim.
\end{itemize}
  
Now, let $q$ be the restriction to $X^{u_t}$ of the canonical projection $\Fl(n) \to G_{b_t,n}$. 
We have to show that $X^{u_{t+1}}=q^{-1}(X^{\ka'_t})$. But since $u_{t}\om_{b_t} \neq \ka'_t$, $u_{t+1}$ is a minimal element of the poset $\{u \ge u_t,\ u\om_{b_t}=\ka'_t\}$, so by Remark~\ref{rem:geomuniquemax} $q^{-1}(X^{\ka'_t})=X^{u_{t+1}}$.~$\square$~\\

\begin{Not}
For every $t \in [0,N-1]$, we set $l_t=j_{b_t}$, that is the largest integer $j$ such that $i_j=b_t$.  
\end{Not}

\begin{Cor}\label{Hypersurfaces}
With the notation of Lemma~\ref{hypersurface_f}, 
\[
Y_\bi^{u_{t+1}}=\{\g \in Z_\bi \mid \g \in Y_\bi^{u_t} \text{\rm\ and\ }p_T(\g)=0\},
\]
where $T=\emptyset*\dots*\ka_t*\dots*\emptyset$ is a tableau of shape $(\bi,0\dots1\dots0)$, the 1 being at position $l_t$.
\end{Cor}

\pr Write $\om=\om_{b_t}$ and $\ka=\ka_t$. Let $\g$ be a gallery
\[
F_\ca\ad{i_1}F_1\ad{i_2}\cdots\ad{b_t}F_{k_t}\ad{} \cdots \ad{} F_r
\]
in $Y_\bi^{u_t}$. This gallery belongs to $Y_\bi^{u_{t+1}}$ if and only if $F_r \in X^{u_{t+1}}$. Since we already know that $F_r \in X^{u_{t}}$, we have
\[
\g \in Y_\bi^{u_{t+1}} \iff p_\ka(F_r)=0 \iff p_\ka(\pi_{b_t}F_r)=0,
\]
where the first equivalence follows from Lemma~\ref{hypersurface_f} and the second from the fact that $\ka$ is of size $b_t$. By definition of $l_t$, no adjacency after $F_{j_t}$ is an $b_t$-adjacency, hence $\pi_{b_t}F_{j_t}=\pi_{b_t}F_{j_t+1}=\dots=\pi_{b_t}F_r$, and therefore,
\[
p_\ka(F_r)=0 \iff p_\ka(F_{j_t})=0 \iff p_T(\g)=0,
\]
where $T=\emptyset*\dots*\ka*\dots*\emptyset$ with $\ka$ in position $l_t$.~$\square$
~\\

\begin{Nots}
We fix an $\ba=a_1\dots a_r$ with $a_i>0$ for every $i$. The associated line bundle $L_{\bi,\ba}$ is very ample, so it induces an embedding of $Z_\bi$ into the projective space $\PP_\ba=\PP(H^0(Z_\bi,L_{\bi,\ba})^*)$ (recall Theorem~\ref{SMT-BS}(3)). We denote by $R_t$ the homogeneous coordinate ring of $Y_\bi^{u_t}$ viewed as a subvariety of $\PP_\ba$.
\end{Nots}

\begin{Rem}
For the rest of this section, if a notion depends on an embedding, such as projective normality, or the homogeneous coordinate ring of a variety, it will be implicitly understood that we work with the line bundle $L_{\bi,\ba}$.   
\end{Rem}

The ring $R_{t+1}$ is a quotient $R_t/I_t$, and we shall determine the ideal $I_t$. We begin by computing the equations of $Y_\bi^{u_{t+1}}$ in an affine open set of $Y_\bi^{u_t}$.

\begin{Def}\label{def_Omega}
We shall define an affine open set $\Omega$ of $Z_\bi$, isomorphic to the affine space $\F^r$. This construction is taken from \cite{Har}.

First, we define inductively a sequence of permutations $(\s_j)$ with $\s_N=w_0$: 
\[
\left\{\begin{array}{ll}
\s_0=e,\\
\s_{j+1}=\s_j*s_{i_{j+1}} & \forall j\ge 0.
\end{array}\right. 
\]
Moreover, we set $v_{j+1}=\s_j^{-1}\s_{j+1} \in \{e,s_{i_{j+1}}\}$. \\

Next, consider the 1-parameter unipotent subgroup $U_\beta$ associated to a root $\beta$, with its standard parametrization $\e_\beta:\F\to U_\beta$ (\emph{i.e.}\ the matrix $\e_\beta(x)$ has 1s on the diagonal, the entry corresponding to $\beta$ equal to $x$, and 0s elsewhere). We also denote by $\al_1,\dots,\al_{n-1}$ the simple roots and by $P_j$ the minimal parabolic subgroup associated to $\al_j$, \emph{i.e.}\ the subgroup generated by $B$ and by $U_{-\alpha_j}$.

We set $\beta_j=v_j(-\al_{i_j})$ and consider the morphism
\[
 \begin{array}{rcl}
  \F^r & \to & P_\bi=P_{i_1}\times \dots \times P_{i_r}\\
 (x_1,\dots,x_r) & \mapsto & (A_1,\dots,A_r)
 \end{array}
\]
with $A_j=\e_{\beta_j}(x_j)v_j$. Set $B_j=A_1\dots A_j$.

Finally, let
\[
\begin{array}{lrcl}
 \phi: & \F^r & \to & Z_\bi\\
 &(x_1,\dots,x_r) & \mapsto & (\g_1,\dots,\g_r)
 \end{array}
\]
for $\g_j=B_j F_\ca$. 

The image of $\phi$ is denoted by $\Omega$: it is an open set in $Z_\bi$, and $\phi: \F^r \to \Omega$ is an isomorphism. 
\end{Def}

\begin{Not}
Let $\ka=k_1\dots k_i$ and $\ta=t_1\dots t_i$ be two columns of the same size. Given a matrix $M$, we denote by $M[\ka,\ta]$ the determinant of the submatrix of $M$ obtained by taking the rows $k_1,\dots,k_i$ and the columns $t_1,\dots,t_i$. Moreover, $M[\ka,[i]]$ is simply denoted by $M[\ka,i]$. 
\end{Not}

\begin{Ex}
We work on Example~\ref{hypersurface_ex}, where $\bi=123212312$, and we set $\ba=111111111$. The sequence $(\s_j)$ is given by
\[
\begin{array}{l}
\begin{array}{@{}ccc}
 \s_0=[1234], &\s_1=[2134], &\s_2=[2314],\\
 \s_3=[2341], &\s_4=[2431], &\s_5=[4231], \\
 \s_6=[4321], & & \\
\end{array} \\
 \s_7=\s_8=\s_9=\s_6.
\end{array}
\]

and the sequence $(v_j)$ is
\[
\begin{array}{l}
\begin{array}{@{}ccc}
v_1=s_1, &v_2=s_2, & v_3=s_3,\\
v_4=s_2, & v_5=s_1, & v_6=s_2, \\
\end{array}\\
v_7=v_8=v_9=e.
\end{array}
\]

Let $T_0=2*23*234*24*4*34*234*4*34$. It can be shown that $\Omega$ is exactly the open set $\{\g \in Z_\bi\mid p_{T_0}(\g)\neq0\}$.\\

Now, direct computations show that 
\[
\phi(x_1,\dots,x_9) \in Y_\bi^{s_1} \iff Q(x_1,\dots,x_9)=x_8(x_1x_6+x_2)+x_1x_5+x_2x_4+x_3=0. 
\]
Since $Y_\bi^{s_1} \cap \Omega$ is irreducible (as an open set of the irreducible $Y_\bi^{s_1}$), the polynomial $Q$ is also irreducible and generates the ideal of the affine vaiety $\phi^{-1}\left(Y_\bi^{s_1} \cap \Omega\right)$. Thus, if $f$ is a linear combination of monomials $p_T$ with $T$ of shape $(\bi,\bm)$ such that $f$ vanishes identically on $Y_\bi^{s_1}$, then $\frac{f}{T_0} \in \F[x_1,\dots,x_9]$ vanishes on $\phi^{-1}\left(Y_\bi^{s_1} \cap \Omega\right)$, hence 
\[
\frac{f}{T_0} \in Q\F[x_1,\dots,x_9]. 
\]
But we know that each coordinate $x_j$ is a quotient $f_j/T_0^{k}$ of degree 0 for an $f_j \in R_0=\F[Z_\bi]$, and also that  
\[
x_8(x_1x_6+x_2)+x_1x_5+x_2x_4+x_3=\frac{p_{T_1}}{p_{T_0}},
\]
where $T_1=2*23*234*24*4*34*234*1*34$. It follows that $f$ is a multiple of $p_{T_1}$, hence $f \in p_1 H^0(Z_\bi,L_{\bi,\ba'})$ where $\ba'=111111101$.

Hence, the ideal $I_1$ of $Y_\bi^{u_1}$ in $R_0$ is $p_1 H^0(Z_\bi,L_{\bi,\ba'})$. We shall generalize this computation below.
\end{Ex}

\begin{Lem}\label{Ubeta}
 For every $j$, 
\[
U_{\beta_1}v_1\dots U_{\beta_j}v_j=U_{\beta_1}U_{\s_1(\beta_2)}\dots U_{\s_{j-1}(\beta_j)}\s_j \subset B\s_j.
\]
\end{Lem}

\pr The equality follows from the formula
\[
\s U_\beta=U_{\s(\beta)}\s,\quad \forall \s \in S_n.
\]
For the inclusion, we proceed by induction over $j$. Since $\beta_1=(i_1,i_1+1)$ and $v_1=s_{i_1}$, $U_{\beta_1}v_1 \subset Bv_1=B\s_1$.

Assume that the property holds for $j\ge 1$, that is $U_{\beta_1}v_1\dots U_{\beta_j}v_j \subset B\s_j$. 

If $\s_js_{i_{j+1}}<\s_j$, then $\s_{j+1}=\s_j$, $v_{j+1}=e$, $\beta_{j+1}=(i_{j+1}+1,i_{j+1})$, and $\s_j(\beta_{j+1})=(\s_j(i_{j+1}+1),\s_j(i_{j+1}))$.
\[
\begin{split}
 \s_js_{i_{j+1}}<\s_j & \iff \s_j(i_{j+1}) > \s_j(i_{j+1}+1)\\
                      & \iff U_{\s_j(\beta_{j+1})} \subset B.
\end{split}
\]
It follows that 
\[
\begin{split}
U_{\beta_1}v_1\dots U_{\beta_j}v_j U_{\beta_{j+1}}v_{j+1} & \subset B\s_jU_{\beta_{j+1}}\\
                                                          & \subset B U_{\s_j(\beta_{j+1})} \s_j\\
                                                          & \subset B\s_{j+1}.
\end{split}
\]

Similarly, if $\s_js_{i_{j+1}}>\s_j$, then $\s_{j+1}=\s_js_{i_{j+1}}$, $v_{j+1}=s_{i_j}$ and $\beta_{j+1}=\al_{i_{j+1}}$. Hence
\[
\begin{split}
U_{\beta_1}v_1\dots U_{\beta_j}v_jU_{\beta_{j+1}}v_{j+1} &\subset B\s_jU_{\bt_{j+1}}s_{i_j}\\
                                                         &\subset BU_{\s_j\bt_{j+1}}\s_js_{i_j}\\
                                                         &\subset B\s_{j+1}.~\square
\end{split}
\]

\begin{Not}
Let $\ka$ be a column of size $i$. We set $O_{\ka}=\{F \in \Fl(n)\mid p_\ka(F)\neq 0\}$. Note that $O_\ka$ is a $T$-stable open subset of $\Fl(n)$.
\end{Not}

\begin{Lem}\label{lem:capO}  
Let $w \in S_n$ and $\ka=w\om_i \in I_{i,n}$. Then $X_w \cap O_\kappa=\coprod_{v\in \E(w,\kappa)}C_v$.
\end{Lem}

\pr Let $q$ be the restriction to $X_w$ of the canonical projection $\Fl(n) \to G_{i,n}$. Then we have $X_w \cap O_\ka =q^{-1}(X_\ka \cap q(O_\ka))$. It is well known that $X_ \ka \cap q(O_\ka)=C_\ka$. But we have the equality 
\[
q^{-1}(C_\ka)=\displaystyle\coprod_{v \in \E(w,\ka)} C_v
\]
(\emph{cf.}\ Remark~\ref{rem:geomuniquemax}).~$\square$

\begin{Prop}\label{prop:local}
There exists a tableau $T_0$ of shape $(\bi,\ba)$ such that 
\[
\Omega=\{\g \in Z_\bi\mid p_{T_0}(\g)\neq 0\}.
 \]
In particular, $\phi$ induces an isomorphism $\phi^*: (R_0)_{(p_{T_0})} \to  \F[x_1,\dots,x_r]$, where $ (R_0)_{(p_{T_0})}$ is the subring of elements of degree~0 in the localized ring $ (R_0)_{p_{T_0}}$, \emph{i.e.}
\[
(R_0)_{(p_{T_0})}=\bigcup_{d\ge 0}\left\{\left. \frac{f}{p_{T_0}^d}\ \right|\ f \in R_0 \textrm{\ is homogeneous of degree $d$}\right\}.
\] 
\end{Prop}

\pr Let $T_0=(\s_1\om_{i_1})^{*a_1}*(\s_2\om_{i_2})^{*a_2}*\dots*(\s_r\om_{i_r})^{*a_r}$. Then 
\[
\{\g \in Z_\bi\mid p_{T_0}(\g)\neq 0\}=\prod_{j=1}^r \proj_j^{-1}\left(O_{\s_j\om_{i_j}}\right).
\]
We know that 
\[
\proj_j(\Omega)=U_{\beta_1}v_1\dots U_{\beta_j}v_jF_\ca.
\]
Thus, by Lemma~\ref{Ubeta}, $\proj_j(\Omega) \subset B\s_jF_\ca=C_{\s_j}$. But if $F \in C_{\s_j}$, then its $i_j$th component $F^{i_j}$ belongs to $C_{\s_j\om_{i_j}}$, so
\[
p_{\s_j\om_{i_j}}(F)=p_{\s_j\om_{i_j}}(F^{i_j})\neq 0.
\] 
This proves the inclusion
\[
\Omega \subset \prod_{j=1}^r \proj_j^{-1}\left(O_{\s_j\om_{i_j}}\right). 
\]

For the opposite inclusion, we proceed by induction over $r$.

If $r=1$, then $\Omega=\{F_\ca\} \times U_{\bt_1}v_1F_\ca$, with $v_1=s_{i_1}$ and $\bt_1=\al_{i_1}$. 
Hence $\Omega=\{F_\ca\} \times C_{s_{i_1}}$. By Lemma~\ref{lem:capO}, $X_{s_{i_1}}\cap O_{s_{i_1}\om_{i_1}}=C_{s_{i_1}}$, so  
\[
\proj_1^{-1}(O_{\s_1\om_{i_1}})=\{F_\ca \ad{i_1} F^1\mid F^1 \in C_{s_{i_1}}\}.
\]
Thus $\Omega=\proj_1^{-1}(O_{\sigma_1\om_{i_1}})$.

Let $r>1$ and assume that the property holds for $r-1$. Let $\g \in \displaystyle\prod_{j=1}^r \proj_j^{-1}\left(O_{\s_j\om_{i_j}}\right).$
By induction, there exist $A_1,\dots,A_{r-1}$ such that $A_j \in U_{\beta_j}v_j$ and $\g_j=A_1\dots A_j F_\ca$ for $j\le r-1$.
Since $\g_{r-1} \ad{i_r} \g_r$, there exists $p \in P_{i_r}$ such that $A_1\dots A_{r-1}pF_\ca=\g_r$. Now, $P_{i_r}F_\ca$ is a 
$T$-stable curve and we have
\[
 P_{i_r}F_\ca=U_{-\al_{i_r}}F_\ca \cup \{s_{i_r}F_\ca\}=U_{\al_{i_r}}s_{i_r}F_\ca \cup \{F_\ca\}.
\]

If $v_r=e$, then $\s_r=\s_{r-1}$ and $\bt_r=-\al_{i_r}$. By Lemma~\ref{Ubeta}, $\g_r \in B\s_{r-1}pF_\ca$. If $pF_\ca=s_{i_r}F_\ca$, then 
$\g_r \in C_{\s_rs_r}$. But by Lemma~\ref{lem:capO}, $\g_r$ belongs to a Schubert cell $C_v$ with $v\om_{i_r}=\s_r\om_{i_r}$. Since 
$\s_r\om_{i_r}\neq \s_rs_r\om_{i_r}$, we have a contradiction. Hence $pF_\ca \in U_{-\al_{i_r}}F_\ca$. So we may choose $p$ in 
$U_{\bt_r}$, which proves that $\g \in \Omega$.

If $v_r=s_{i_r}$, we prove similarly that $pF_\ca \neq F_\ca$, so we may choose $p$ in $U_{\al_{i_r}}s_{i_r}=U_{\bt_r}v_r$, thus 
$\g \in \Omega$.~$\square$

\begin{Rem}
Consider an arbitrary tableau $T$ of shape $(\bi,\ba)$. Then we may compute $\phi^*\left(\frac{p_T}{p_{T_0}}\right)$ in the following way. Write 
\[
T=\ka_{1,1}*\dots*\ka_{1,a_1}*\dots*\ka_{r,1}*\dots*\ka_{r,a_r},
\] 
then 
\[
 \phi^*\left(\frac{p_T}{p_{T_0}}\right)(x_1,\dots,x_r)=B_1[\ka_{1,1},i_1]\dots B_1[\ka_{1,a_1},i_1]\dots B_r[\ka_{r,1},i_r]\dots B_r[\ka_{r,a_r},i_r].
\]
\end{Rem}
\vspace*{2mm}

\begin{Prop} The entries of $B_r$ below the antidiagonal are zero. Denote by $Q_{i,j} \in \F[x_1,\dots,x_r]$ the coefficients of $B_r$ above the antidiagonal: 
\[
 B_r=\left(\begin{array}{cccccc}
            Q_{1,1} & Q_{1,2} & \cdots &\cdots & Q_{1,n-1} & 1\\
            \vdots & \vdots & & &1 & 0\\
            \vdots & \vdots & & \iddots&0 & 0\\
            \vdots & \vdots &\iddots&\iddots  & \vdots & \vdots\\
            Q_{n-1,1} & 1 &\iddots  &&\vdots & \vdots\\
            1 & 0 & \cdots & \cdots& 0& 0
           \end{array}\right).
\]

Then the polynomials $Q_{i,j}$ all have distinct nonzero linear parts.
\end{Prop}

\pr  From Lemma~\ref{Ubeta}, we have 
\[
B_r \in U_{\beta_1}v_1\dots U_{\beta_r}v_rF_\ca\subset B \s_r=Bw_0,
\]
so only the coefficients of $B_r$ above the antidiagonal are nonzero.

We shall obtain the linear part of $Q_{i,j}$ by differentiating $B_r$. From the expression $B_r=\e_{\beta_1}(x_1)v_1\dots\e_{\beta_r}(x_r)v_r$, we see that 
\[
\frac{\partial B_r}{\partial x_k}(0)=E_{\s_{k-1}(\beta_k)}w_0. 
\]

This proves that the linear part of $Q_{i,j}$ is the sum of all $x_k$ such that $\s_{k-1}\beta_k=(i,n+1-j)$. Hence all these linear parts are distinct, provided they are nonzero, which amounts to the fact that
\[
U_{\beta_1}U_{\s_1\beta_2}\dots U_{\s_{r-1}\beta_r}=U
\] 
(where $U$ is the unipotent part of $B$). Since the stabilizer of $w_0F_\ca$ in $U$ is trivial, this is equivalent to
\[
U_{\beta_1}U_{\s_1(\beta_2)}\dots U_{\s_{r-1}(\beta_r)}w_0F_\ca=Uw_0F_\ca,
\]
or, using Lemma~\ref{Ubeta}, equivalent to
\[
U_{\beta_1}v_1\dots U_{\beta_r}v_rF_\ca=Uw_0F_\ca,
\]
which indeed holds because $\proj_r(\Omega)=C_{w_0}$.~$\square$
\vspace*{2mm}

\begin{Prop}~\label{local} 
The ideal of $Y_\bi^{u_{t+1}}\cap \Omega$ in the coordinate ring of $Y_\bi^{u_t} \cap \Omega$ is generated by $Q_{\ka_{t,1},b_t}$,
where $\ka_{t,1}$ is the first entry of the column $\ka_t$. Moreover, 
\[Q_{\ka_{t,1},b_t}=\phi^*\left(\frac{p_{T_t}}{p_{T_0}}\right)\]
where $T_t$ is the tableau obtained from $T_0$ by replacing its last column of size $b_t$ by $\ka_t$.

Moreover, the varieties $Y_\bi^{u_t} \cap \Omega$ are isomorphic to affine spaces.
\end{Prop}

\pr We already know by Corollary~\ref{Hypersurfaces} that, for a gallery $\g \in Y_\bi^{u_t}$,
\[
 g \in Y_\bi^{u_{t+1}} \text{\ if and only if\ } p_{\ka_t}(F_{l_t})=0.
\]

In $\Omega$, this corresponds to the vanishing of $B_{l_t}[\s_{l_t}\om_{b_t},b_t]$. Since $i_j\neq b_t$ for $j>l_t$, the spaces generated by the first $b_t$ columns of $B_{l_t}$ and of $B_r$ are equal, hence
\[
 B_{l_t}[\s_{l_t}\om_{b_t},b_t]=0 \iff  B_r[\s_{l_t}\om_{b_t},b_t]=0.
\]
A straightforward computation shows that $B_r[\s_{l_t}\om_{b_t},b_t]=\pm Q_{\ka_{t,1},b_t}$. 

To prove that $\phi^*\left(\frac{p_{T_t}}{p_{T_0}}\right)=\pm Q_{\ka_{t,1},b_t}$, note that 
\[
\phi^*\left(\frac{p_{T_t}}{p_{T_0}}\right)=B_1[\s_1\om_{i_1},i_1]\dots B_{j_t}[\s_{j_t}\om_{b_t},b_t]\dots B_r[\s_r\om_{i_r},i_r]. 
\]
Now, by Lemma~\ref{Ubeta}, $B_j=b_j\s_j$ for some $b_j \in B$. So, for $j\neq l_t$, 
\[
B_j[\s_j\om_{i_j},i_j]=\pm b_j[\s_j\om_{i_j},\s_j\om_{i_j}]=\pm 1. 
\]

Hence $\phi(x_1,\dots,x_r) \in Y_\bi^{u_{t+1}}\cap \Omega$ if and only if $Q_{\ka_{t,1},b_t}(x_1,\dots,x_r)=0$. But this polynomial is of the form $x_{p_t}-Q'$ for some variable $x_{p_t} \in \{x_1,\dots,x_r\}$, so we may substitute $Q'$ for $x_{p_t}$ in the coordinate ring of $Y_\bi^{u_t} \cap \Omega$ to obtain the coordinate ring of $Y_\bi^{u_{t+1}}\cap \Omega$. Thus, by induction over $t$, we see that the coordinate ring of $Y_\bi^{u_{t+1}} \cap \Omega$ is isomorphic to $\F[x_i\mid i\neq p_0,\dots,p_t]$. In particular, this ring is a Unique Factorization Domain. Therefore, the irreducible polynomial $Q_{\ka_{t,1},b_t}$ generates the ideal of $Y_\bi^{u_{t+1}}\cap \Omega$ in the coordinate ring of $Y_\bi^{u_t}\cap \Omega.~\square$

\begin{Nots}
We set $\ba'_t=0\ldots-1\ldots0$, the $-1$ again being at position $l_t$. Let $S_t$ be the $R_t$-graded module associated to the coherent sheaf $L_{\bi,\ba'_t}$, that is,
\[
S_t=\bigoplus_{d=0}^{+\infty} H^0(Y_\bi^{u_t},L_{\bi,d\ba+\ba'_t}).
\]
\end{Nots}

\begin{Hyp}\label{hyp:regular}
We assume that $\ba+\displaystyle\sum_{t=0}^{N-1} \ba'_t$ is regular.
\end{Hyp}

\begin{Cor}~\label{exactsequence}
Denote by $\OO_{Y_\bi^{u_t}}$ the structural sheaf of $Y_\bi^{u_t}$ and assume that $Y_\bi^{u_t}$ is projectively normal. Then the sequence of $R_t$-modules
\[
(*) \qquad 0 \to S_t \to R_t \to R_{t+1} \to 0
\]
is exact, where the first map is the multiplication by $p_{\ka_t}$ and the second is the natural projection.

The exact sequence $(*)$ induces an exact sequence of sheaves of $\OO_{Y_\bi^{u_t}}$-modules 
\[
 0 \to L_{\bi,\bm+\ba'_t} \to L_{\bi,\bm} \to (L_{\bi,\bm})_{|Y_\bi^{u_{t+1}}} \to 0
\] 
and a long exact sequence in cohomology
 \[
(**) \qquad \xymatrix@C=5mm{
            0\ar[r] & H^0(Y_\bi^{u_t},L_{\bi,\bm+\ba'_t}) \ar[r] & H^0(Y_\bi^{u_t},L_{\bi,\bm}) \ar[r] & H^0(Y_\bi^{u_{t+1}},L_{\bi,\bm}) \ar `r[d] `_l[l] `^d[dlll] `^r[dll] [dll] \\
                    & H^1(Y_\bi^{u_t},L_{\bi,\bm+\ba'_t}) \ar[r] & \cdots                           &                                                                         \\
                    &                                   & \cdots \ar[r] & H^{i-1}(Y_\bi^{u_{t+1}},L_{\bi,\bm}) \ar `r[d] `_l[l] `^d[dlll] `^r[dll] [dll] \\
                    & H^i(Y_\bi^{u_t},L_{\bi,\bm+\ba'_t}) \ar[r] & H^i(Y_\bi^{u_t},L_{\bi,\bm}) \ar[r] & H^i(Y_\bi^{u_{t+1}},L_{\bi,\bm}) \ar `r[d] `_l[l] `^d[dlll] `^r[dll] [dll] \\
                    & H^{i+1}(Y_\bi^{u_t},L_{\bi,\bm+\ba'_t}) \ar[r] & \cdots                            &
}
\]
\end{Cor}

\pr Since $Y_{\bi}^{u_t}$ is projectively normal, we know that 
\[
R_t=\bigoplus_{d=0}^{+\infty} H^0(Y_\bi^{u_t},L_{\bi,d\ba}), 
\]
hence the sequence 
\[
\xymatrix{
 0 \ar[r] & S_t \ar[r]^{\mu}&  R_t \ar[r]^{q} &  R_{t+1} \ar[r] & 0
}\]
is well defined. Moreover, we already know by Corollary~\ref{Hypersurfaces} that $q \circ \mu =0$. 

Let $f$ be a homogeneous element of degree $d$ in $R_t$, and suppose that $q(f)=0$, that is, $f$ vanishes identically on $Y_\bi^{u_{t+1}}$. Then $\frac{f}{p_{T_0}^d}$ vanishes identically on $Y_\bi^{u_{t+1}}\cap \Omega$, hence $\phi^*\left(\frac{f}{p_{T_0}^d}\right)\in \F[x_1,\dots,x_r]$ is a multiple of $Q_{\ka_{t,1},b_t}=\phi^*\left(\frac{p_{T_t}}{p_{T_0}}\right)$. It follows that $f$ is a multiple of $p_{T_t}$, hence $f \in p_{\ka_t}S_t=\mu(S_t)$.

If we consider the coherent sheaves associated to these $R_t$-modules and tensor them by $L_{\bi,\bm}$, then we get the exact sequence of sheaves of $\OO_{Y_\bi^{u_t}}$-modules
\[
 0 \to L_{\bi,\bm+\ba'_t} \to L_{\bi,\bm} \to (L_{\bi,\bm})_{|Y_\bi^{u_{t+1}}} \to 0,
\]
which gives the long exact sequence $(**)$.~$\square$
\vspace*{2mm}

\begin{Def}
A line bundle $L_{\bi,\bm}$ (or just $\bm$) is said to be \emph{nonnegative} if $m_i \ge 0$ for every $i$.
\end{Def}

\begin{The}~\label{cohomology} For every $t\le N$,
\begin{enumerate}[$(1)_t$]
 \item the variety $Y_\bi^{u_t}$ is projectively normal relatively to $L_{\bi,\ba}$, under Assumption~\ref{hyp:regular},
 \item for every $i>0$, and every $\bm$ such that $\bm^{[t]}=\bm+\ba'_0+\dots+\ba'_{t-1}$ is nonnegative, $H^i(Y_\bi^{u_t},L_{\bi,\bm})=0$. In particular, $H^i(Y_\bi^{u_t},L_{\bi,d\ba})=0$ for every $d \in \Z_{\ge 0}$.
\end{enumerate}
\end{The}

\pr We proceed by induction over $t$. 

For $t=0$, $Y_\bi^{u_0}=Z_\bi$. By Theorem~\ref{SMT-BS}, $Z_\bi$ is projectively normal, and we have $H^i(Z_\bi,L_{\bi,\bm})=0$ for $i>0$. 

Assume that the theorem is true for a $t\ge 0$.\\

We shall prove that $(1)_{t+1}$ is true. By induction, $Y_\bi^{u_t}$ is projectively normal, so the sequence
\[
 0 \to H^0(Y_\bi^{u_t},L_{\bi,\bm+\ba'_t}) \to H^0(Y_\bi^{u_t},L_{\bi,\bm}) \to H^0(Y_\bi^{u_{t+1}},L_{\bi,\bm}) \to H^1(Y_\bi^{u_t},L_{\bi,\bm+\ba'_t}) 
\]
is exact. Moreover, for $d \in \Z_{\ge 0}$, by $(2)_t$ we have $H^1(Y_\bi^{u_t},L_{\bi,d\ba})=0$, hence an exact sequence
\[
 0 \to S_t \to R_t \to \bigoplus_{d\ge 0} H^0(Y_\bi^{u_{t+1}},L_{\bi,d\ba}) \to 0.
\]
Since the sequence 
\[
 0 \to S_t \to R_t \to R_{t+1} \to 0
\]
is also exact, we have $R_{t+1}=\displaystyle\bigoplus_{d\ge 0} H^0(Y_\bi^{u_{t+1}},L_{\bi,d\ba})$, that is, $Y_\bi^{u_{t+1}}$ is projectively normal.\\

We now prove that $(2)_{t+1}$ is true. Let $\bm$ be such that $\bm^{[t+1]}$ is non negative. Note that $\bm^{[t]}$ is also non negative. Since $Y_\bi^{u_t}$ is projectively normal, we have the exact sequence
\[
   H^i(Y_\bi^{u_t},L_{\bi,\bm}) \to H^i(Y_\bi^{u_{t+1}},L_{\bi,\bm}) \to H^{i+1}(Y_\bi^{u_t},L_{\bi,\bm+\ba'_t}).                                                                                                                          
\]
Now, by $(2)_t$,
\[
 \begin{array}{l@{\ \implies\ }l}
  \bm^{[t]} \textrm{\ non negative\ } & H^i(Y_\bi^{u_t},L_{\bi,\bm})=0,\\
  \bm^{[t+1]} \textrm{\ non negative\ } & H^{i+1}(Y_\bi^{u_t},L_{\bi,\bm+\ba'_t})=0.
 \end{array}
\]
Thus, $H^i(Y_\bi^{u_{t+1}},L_{\bi,\bm})=0$.~$\square$

\begin{Cor}
For every $t\le N-1$, if $\bm^{[t+1]}$ is non negative, then the restriction map $H^0(Y_\bi^{u_t},L_{\bi,\bm}) \to H^0(Y_\bi^{u_{t+1}},L_{\bi,\bm})$ is surjective.
\end{Cor}

\pr Since $Y_\bi^{u_{t+1}}$ is projectively normal, we have the exact sequence
\begin{multline}
 0 \to H^0(Y_\bi^{u_{t+1}},L_{\bi,\bm+\ba'_{t+1}}) \to H^0(Y_\bi^{u_{t+1}},L_{\bi,\bm})
\\
 \to H^0(Y_\bi^{u_{t+2}},L_{\bi,\bm}) \to H^1(Y_\bi^{u_{t+1}},L_{\bi,\bm+\ba'_{t+1}}).
\end{multline}
Since $\bm^{[t+1]}$ is regular, we have by $(2)_{t+1}$ that $H^1(Y_\bi^{u_{t+1}},L_{\bi,\bm+\ba'_{t+1}})=0$, hence the restriction map $H^0(Y_\bi^{u_{t+1}},L_{\bi,\bm}) \to H^0(Y_\bi^{u_{t+2}},L_{\bi,\bm})$ is surjective.~$\square$

\begin{Cor}\label{basis-regular}
If $\bm^{[N]}=\bm+\displaystyle\sum_{t=0}^{N-1} \ba'_t$ is non negative, then a basis of $H^0(\RR_{\bi},L_{\bi,\bm})$ is given by the $w_0$-standard monomials of shape $(\bi,\bm)$.
\end{Cor}

\pr Since the restriction $H^0(Z_\bi,L_{\bi,\bm}) \to H^0(\RR_{\bi},L_{\bi,\bm})$ is surjective, the standard monomials $p_T$ that do not vanish identically on $\RR_\bi$ form a generating set. By Theorem~\ref{non_vanishing}, these monomials are exactly the $w_0$-standard monomials. By Theorem~\ref{Linear independence}, these monomials are linearly independent.~$\square$

\begin{Prop}
Let $p_T$ be a standard monomial of shape $(\bi,\bm)$, with $\bm$ arbitrary. Then $p_T$ decomposes as a linear combination of $w_0$-standard monomials on $\RR_\bi$.
\end{Prop}

\pr With the notation of Theorem~\ref{non_vanishing}, the result is true if $\bm^{[N]}$ is non negative. If this is not the case, then we set $\bb=b_1\dots b_r$ with
\[
 b_j=\left\{\begin{array}{ll}
             N & \text{\ if $j=l_t$ for some $t$},\\
             0 & \textrm{otherwise},
            \end{array}\right.
\]
and $\bm'=\bm+\bb$. Now, $\bm'$ satisfies the assumption of Theorem~\ref{non_vanishing}. We multiply $p_T$ by $p_{\ka_1}^N\dots p_{\ka_{n-1}}^N$ to obtain a new monomial $p_T'$, where $\ka_j=w_0\om_{j}$. So $p_T'$ is of shape $(\bi,\bm')$ and does not vanish identically on $\RR_\bi$. Now, $p_T'$ may be non standard, so we decompose it as a linear combination of $w_0$-standard monomials of shape $(\bi,\bm')$ on $\RR_\bi$, thanks to Corollary~\ref{basis-regular}:
\[
 p_T(p_{\ka_1}\dots p_{\ka_{n-1}})^N = p_{T'}=\sum_{T''} a_{T''}\,p_{T''}.
\]

Since a $w_0$-standard monomial does not vanish identically on $\RR_\bi$, the columns $\ka$ that are in position $j_k$ in a tableau $T''$ are maximal, \emph{i.e.} equal to $w_0\om_{k}$. Hence we may factor this linear combination by $(p_{\ka_1}\dots p_{\ka_{n-1}})^N$, so that $p_T$ is a linear combination of $w_0$-standard monomials.~$\square$

\begin{Cor}~
\begin{enumerate}
\item For an arbitrary $\bm$, a basis of $H^0(\RR_{\bi},L_{\bi,\bm})$ is given by the $w_0$-standard monomials of shape $(\bi,\bm)$.
\item For $\bm$ regular, the variety $\RR_\bi$ is projectively normal relatively to the line bundle $L_{\bi,\bm}$.
\end{enumerate}
\end{Cor}

\pr The first point is immediate. For the second point, assume that $\bm$ is regular. We already know by Corollary~\ref{cor-basis_coord_ring} that the $w_0$-standard monomials form a basis of the homogeneous coordinate ring $R_N$ of $\RR_\bi$ in $\PP(H^0(Z_\bi,L_{\bi,\bm})^*)$. Hence $R_N$ and $\bigoplus_d H^0(\RR_\bi,L_{\bi,\bm}^{\otimes d})$ are equal.~$\square$

\begin{Rem}~
\begin{enumerate}
\item In the regular case ($m_i>0$ for every $i$), the basis given by standard monomials is compatible with $\RR_\bi$: this is no longer the case if $\bm$ is not regular, see Remark~\ref{counter_example}.
\item By Theorem~\ref{cohomology}, we know that if $\bm^{[N]}$ is non negative, then $H^i(\RR_\bi,L_{\bi,\bm})=0$ for all $i>0$. I do not know if this vanishing result is still valid for an arbitrary $\bm$.
\item Similarly, if $\bm$ is not regular, I do not know if $\RR_\bi$ is projectively normal.
\end{enumerate}
\end{Rem}

\begin{Rem}
Since the notion of standard tableau for Bott-Samelson varieties is defined in types other than~A (see \cite{LLM2}), one may ask whether the results of this paper extend to this more general setting. The Demazure product is defined in arbitrary type, which could again lead to good properties for $w_0$-standard tableaux. However, arguments used in the proofs of Theorem~\ref{Linear independence} and Theorem~\ref{non_vanishing} are specific to type~A.
\end{Rem}

\end{document}